\documentclass[10pt]{article}
\usepackage{amsmath}
\usepackage{amsthm,fullpage}
\usepackage{amssymb}
\usepackage{epsfig}
\usepackage[latin1]{inputenc}

\newtheorem{thm}{Theorem}[section]
\newtheorem{prop}[thm]{Proposition}
\newtheorem{cor}[thm]{Corollary}
\newtheorem{lemma}[thm]{Lemma}

\theoremstyle{definition} 
\newtheorem{Def}[thm]{Definition}

\newcommand{\Pa}{\textrm{\textbf{p}}}

\newcommand{\DD}{\textbf{D}}
\newcommand{\NC}{\textbf{NC}}
\newcommand{\BB}{\textbf{B}}
\newcommand{\LS}{\mathcal{L}^S}
\newcommand{\LT}{\mathcal{L}^T}
\newcommand{\LK}{\mathcal{L}^K}
\newcommand{\mW}{\mathcal{W}}
\newcommand{\lS}{\leq_S}
\newcommand{\lT}{\leq_T}
\newcommand{\lK}{\leq_K}

\newcommand{\BN}{\mathcal{\textbf{N}}}
\newcommand{\BS}{\mathcal{\textbf{S}}}
\newcommand{\PR}{P\!R}

\newcommand{\commentaire}[1]{}

\newcommand{\Ref}[1]{(\ref{#1})}
\newcommand{\dem}{\noindent \textbf{Proof: }}
\newcommand{\demof}[1]{\noindent \textbf{Proof of~{#1}:}}
\newcommand{\findem}{\vspace{-.55cm} \begin{flushright} $\square~$ \end{flushright} \vspace{.2cm} }
\newcommand{\findembis}[1]{\hspace{#1} \vspace{.1cm} \raisebox{-.15cm}{$\square~$}\\}
\newcommand{\bbar}[1]{\overline{#1}}

\newcommand{\titre}[1]{\noindent \textbf{#1}}
\newcommand{\ite}{\noindent $\bullet~$}
\newcommand{\iten}{\noindent -~}

\newcommand{\al}{\alpha} 
\newcommand{\be}{\beta} 
\newcommand{\under}[1]{\raisebox{.05cm}{~\underline{#1}~}}
\newcommand{\underder}[1]{\raisebox{.1cm}{~\underline{\underline{#1}}~}}

\newcommand{\tth}{\textrm{th}}

\title{Catalan's intervals and realizers of triangulations}

\author{Olivier Bernardi and Nicolas Bonichon}

\begin{document}




\maketitle
\begin{abstract}
The Stanley lattice, Tamari lattice and Kreweras lattice are three remarkable orders defined on the set of Catalan objects of a given size. These lattices are ordered by inclusion: the Stanley lattice is an extension of the Tamari lattice which is an extension of the Kreweras lattice. The Stanley order can be defined on the set of Dyck paths of size $n$ as the relation of \emph{being above}. Hence, intervals in the Stanley lattice are pairs of non-crossing Dyck paths. In a former article, the second author defined a bijection $\Phi$ between pairs of non-crossing Dyck paths and the realizers of triangulations (or Schnyder woods). We give a simpler description of the bijection $\Phi$. Then, we study the restriction of $\Phi$ to Tamari's and Kreweras' intervals. We prove that $\Phi$ induces a bijection between Tamari intervals and minimal realizers. This gives a bijection between Tamari intervals and triangulations. We also prove that  $\Phi$ induces a bijection between Kreweras intervals and the (unique) realizers of stack triangulations. Thus, $\Phi$ induces a bijection between Kreweras intervals and stack triangulations which are known to be in bijection with ternary trees. \vspace{-.2cm}
\end{abstract}

\section{Introduction}
\begin{figure}[htb!]
\begin{center}
\input{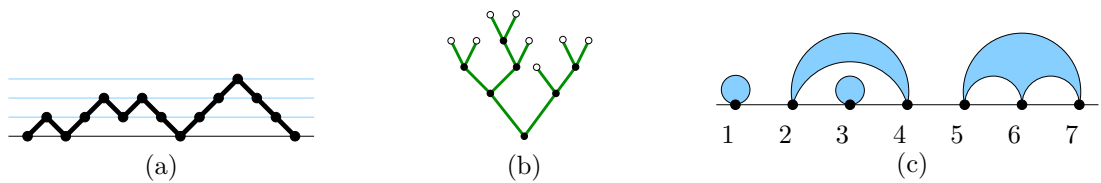}
\caption{(a) A Dyck path. (b) A binary tree. (c) A non-crossing partition.}\label{fig:Dyck-exp}
\end{center}
\end{figure}

\begin{figure}[htb!]
\begin{center}
\input{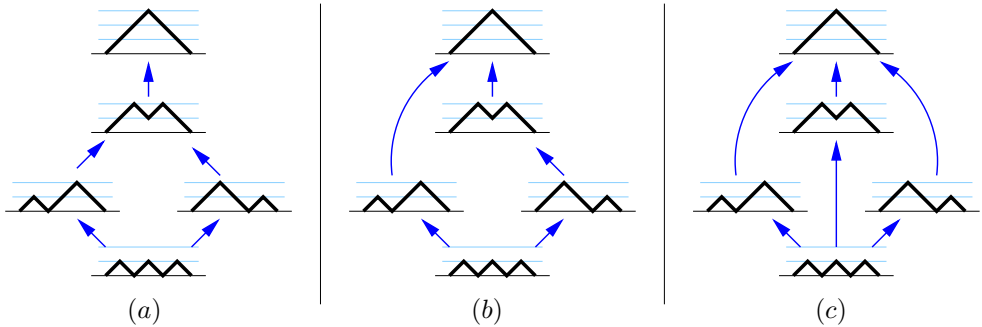}
\caption{Hasse diagrams of the Catalan lattices on the set $\DD_3$ of Dyck paths:  $(a)$ Stanley lattice, $(b)$ Tamari lattice, $(c)$ Kreweras lattice.}\label{fig:Catalan-lattices}
\end{center}
\end{figure}

A \emph{Dyck path} is a lattice path made of  $+1$ and $-1$ steps that starts from 0, remains non-negative and ends at 0. It is often convenient to represent a Dyck path by a sequence of  North-East and South-East steps as is done in Figure~\ref{fig:Dyck-exp}~(a).
The set $\DD_n$ of Dyck paths of length $2n$ can be ordered by the relation $P\lS Q$ if $P$ stays below $Q$. This partial order is in fact a distributive lattice on $\DD_n$ known as the \emph{Stanley lattice}. The Hasse diagram of the Stanley lattice on $\DD_3$ is represented in  Figure~\ref{fig:Catalan-lattices} $(a)$.

It is well known that the Dyck paths of length $2n$ are counted by the $n^{\textrm{th}}$ \emph{Catalan number}  $C_n=\frac{1}{n+1}{2n \choose n}$. The Catalan sequence is a pervasive guest in enumerative combinatorics. Indeed, beside Dyck paths, this sequence enumerates the binary trees, the plane trees, the non-crossing partitions and over 60 other fundamental combinatorial structures \cite[Ex. 6.19]{Stanley:volume2}. These different incarnations of the Catalan family gave rise to several lattices beside Stanley's. The \emph{Tamari lattice} appears naturally in the study of binary trees where the covering relation corresponds to right rotation. This lattice is actively studied due to its link with the associahedron (Stasheff polytope). Indeed, the Hasse diagram of the Tamari lattice is the 1-skeleton of the associahedron. The \emph{Kreweras lattice} appears naturally in the setting of non-crossing partitions.  In the seminal paper \cite{Kreweras:non-crossing}, Kreweras proved that the refinement order on non-crossing partitions defines a lattice. Kreweras lattice appears to support a great deal of mathematics that reach far beyond enumerative combinatorics \cite{McCammond:non-crossing-partitions,Simion:non-crossing-partition}. Using suitable bijection between Dyck paths, binary trees, non-crossing partitions and plane trees,  the three \emph{Catalan lattices} can be defined on the set of plane trees of size $n$ in such way that the Stanley lattice $\LS_n$ is an extension of the Tamari lattice $\LT_n$ which in turn is an extension of the Kreweras lattice $\LK_n$ (see \cite[Ex. 7.2.1.6 - 26, 27 and 28]{Knuth:generating-all-trees}). In this paper, we shall find convenient to embed the three Catalan lattices on the set $\DD_n$ of Dyck paths. The Hasse diagram of the Catalan lattices on $\DD_3$ is represented in Figure~\ref{fig:Catalan-lattices}. \\

There are closed formulas for the number of \emph{intervals}  (i.e. pairs of comparable elements) in each of the Catalan lattices. The intervals of the Stanley lattice are the pairs of non-crossing Dyck paths and  the number $|\LS_n|$ of such pairs can be calculated using the lattice path determinant formula of Lindström-Gessel-Viennot \cite{Gessel-Viennot:lattice-determinant}. It is shown in \cite{Desainte-Catherine:enumeration-Young-Tableaux} that 
\begin{eqnarray}\label{eq:Stanley-intervals}
|\LS_n|~=~C_{n+2}C_{n} - C_{n+1}^2~=~\frac{6 (2n)! (2n+2)!}{n! (n+1)! (n+2)! (n+3)!}.  
\end{eqnarray}
The intervals of the Tamari lattice were recently enumerated by Chapoton \cite{Chapoton:Tamari-intervals} using a generating function approach. It was proved that the number  of intervals in the Tamari lattice is  
\begin{eqnarray}\label{eq:Tamari-intervals}
\displaystyle |\LT_n|~=~\frac{2(4n+1)!}{(n+1)!(3n+2)!}.
\end{eqnarray}
Chapoton also noticed that~\Ref{eq:Tamari-intervals} is the number of triangulations (i.e. maximal planar graphs)  and asked for an explanation.  
The number $|\LK_n|$ of intervals of the Kreweras Lattice has an even simpler formula. In \cite{Kreweras:non-crossing}, Kreweras proved by a recursive method that 
\begin{eqnarray}\label{eq:Kreweras-intervals}
\displaystyle |\LK_n|~=~\frac{1}{2n+1}{3n \choose n}.
\end{eqnarray}
This is also the number of ternary trees and a bijection was exhibited in \cite{Edelman:non-crossing-chains-trees}. 
\\

In \cite{Bonichon:realizers}, the second author defined a bijection $\Phi$ between  the pairs of non-crossing Dyck paths (equivalently, Stanley's intervals) and the \emph{realizers} (or \emph{Schnyder woods}) of triangulations. The main purpose of this article is to study the restriction of the bijection $\Phi$ to the Tamari intervals and to the Kreweras intervals. We first give an alternative, simpler, description of the bijection $\Phi$. Then, we prove that the bijection $\Phi$ induces a bijection between the intervals of the Tamari lattice and the realizers which are \emph{minimal}. Since every triangulation has a unique \emph{minimal} realizer, we obtain a bijection between Tamari intervals and triangulations. As a corollary, we obtain a bijective proof of Formula~\Ref{eq:Tamari-intervals} thereby answering the question of Chapoton. Turning to the Kreweras lattice, we prove that the mapping $\Phi$ induces a bijection between Kreweras intervals and the realizers which are both  \emph{minimal} and \emph{maximal}. We then characterize the triangulations having a realizer which is both minimal and maximal and prove that these triangulations are in bijection with ternary trees. This gives a new bijective proof of Formula~\Ref{eq:Kreweras-intervals}.\\


The outline of this paper is as follows. In Section~\ref{section:Catalan-lattice}, we review our notations about Dyck paths and characterize the covering relations for the Stanley, Tamari and Kreweras lattices in terms of Dyck paths. In Section~\ref{section:bijection}, we recall the definitions about triangulations and realizers. We then give an alternative description of the bijection $\Phi$ defined in \cite{Bonichon:realizers} between  pairs of non-crossing Dyck paths and the realizers. In Section~\ref{section:Tamari}, we study the restriction of $\Phi$ to the Tamari intervals. Lastly, in Section~\ref{section:Kreweras} we study the restriction of $\Phi$ to the Kreweras intervals.

\section{Catalan lattices}\label{section:Catalan-lattice}
\titre{Dyck paths.} A \emph{Dyck path} is a lattice path made of steps $N=+1$ and $S=-1$ that starts from 0, remains non-negative and ends at 0. A Dyck path is said to be \emph{prime} if it remains positive between its start and end. The \emph{size} of a path is half its length and the set of Dyck paths of size $n$ is denoted by $\DD_n$. \\

\begin{figure}[htb!]
\begin{center}
\input{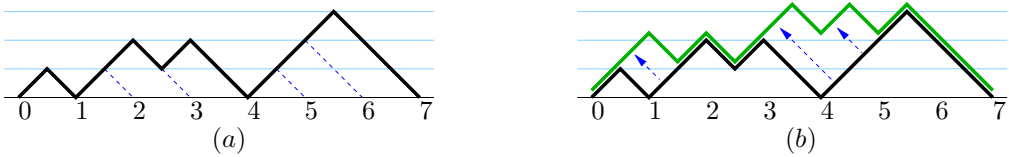}
\caption{(a) Exceedence of a Dyck path. (b) Differences between two Dyck paths.}\label{fig:notations-paths}
\end{center}
\end{figure}

Let $P$ be a Dyck path of size $n$. Since $P$ begins by an $N$ step and has $n$ $N$ steps, it can be written as $P=NS^{\al_1}NS^{\al_2}\ldots NS^{\al_n}$. We call $i^\tth$ \emph{descent} the subsequence $S^{\al_i}$ of $P$.  For $i=0,1,\ldots,n$ we call $i^\tth$ \emph{exceedence} and denote by $e_i(P)$ the height of the path $P$ after the $i^\tth$ descent, that is, $ e_i(P)=i-\sum_{j\leq i}\al_j$. For instance, the Dyck path represented in Figure~\ref{fig:notations-paths}~(a) is $P=NS¹NS^0NS¹NS²NS^0NS^0NS³$ and $e_0(P)=0$, $e_1(P)=0$, $e_2(P)=1$, $e_3(P)=1$,  $e_4(P)=0$,  $e_5(P)=1$, $e_6(P)=2$ and $e_7(P)=0$. If $P,Q$ are two Dyck paths of size $n$, we denote $\delta_i(P,Q)=e_i(Q)-e_i(P)$ and $\Delta(P,Q)=\sum_{i=1}^n \delta_i(P,Q)$. For instance, if $P$ and $Q$ are respectively the lower and upper paths in Figure~\ref{fig:notations-paths}~(b), the values $\delta_i(P,Q)$ are zero except for $\delta_1(P,Q)=1$, $\delta_4(P,Q)=2$ and $\delta_5(P,Q)=1$.

For $0\leq i\leq j\leq n$, we write  $i \under{P} j$ (resp. $i\underder{P}j$) if $e_i(P)\geq e_j(P)$ and $e_i(P)\leq e_k(P)$ (resp. $e_i(P)< e_k(P)$) for all $i<k<j$. In other words, $i \under{P} j$ (resp. $i\underder{P}j$) means that the subpath  $NS^{\al_{i+1}}NS^{\al_{i+2}}\ldots NS^{\al_j}$ is a Dyck path (resp. prime Dyck path) followed by $e_i(P)-e_j(P)$ $S$ steps. For instance, for the Dyck path $P$ of Figure~\ref{fig:notations-paths}~(a), we have $0\under{P} 4$, $~1\underder{P} 4$ and $2\under{P} 4$ (and many other relations). \\

We will now define the Stanley, Tamari and Kreweras lattices in terms of Dyck paths. More precisely, we will characterize the covering relation of each lattice in terms of Dyck paths and show that our definitions respects the known hierarchy between the three lattices (the Stanley lattice is a refinement of the Tamari lattice which is refinement of the Kreweras Lattice; see \cite[Ex. 7.2.1.6 - 26, 27 and 28]{Knuth:generating-all-trees}).\\

\begin{figure}[ht]
\begin{center}
\input{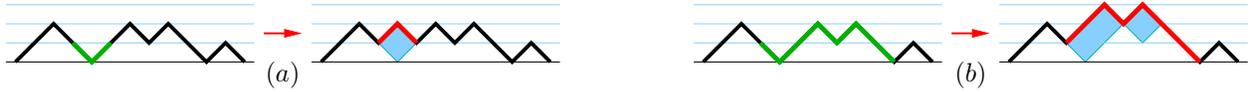}
\caption{Covering relations in (a) Stanley lattice, (b) Tamari lattice.}\label{fig:covering}
\end{center}
\end{figure}

\titre{Stanley lattice.} Let $P=NS^{\al_1}\ldots NS^{\al_n}$ and $Q=NS^{\be_1}\ldots NS^{\be_n}$ be two Dyck paths of size $n$. We denote by $P \lS Q$ if the path $P$ stays below the path $Q$. Equivalently, $e_i(P)\leq e_i(Q)$ for all $1\leq i\leq n$. The relation $\lS$ defines the \emph{Stanley lattice} $\LS_n$ on the set $\DD_n$. Clearly the path $P$ is covered by the path $Q$ in the Stanley lattice if $Q$ is obtained from $P$ by replacing a subpath $SN$ by $NS$. Equivalently, there is an index $1\leq i\leq n$ such that $\be_i=\al_i-1$, $\be_{i+1}=\al_{i+1}+1$ and $\be_k=\al_k$ for all $k\neq i,i+1$. 
The covering relation of the Stanley lattice is represented in Figure~\ref{fig:covering}~(a) and the Hasse Diagram of $\LS_3$ is represented in Figure~\ref{fig:Catalan-lattices}~(a).\\ 



\begin{figure}[ht]
\begin{center}
\input{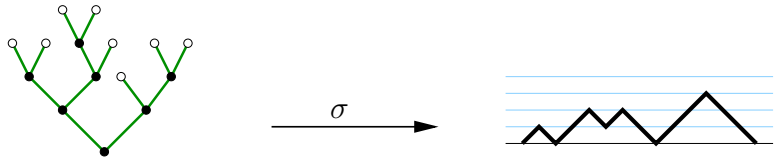}
\caption{The binary tree $(((\circ,\circ),((\circ,\circ),\circ)),(\circ,(\circ,\circ)))$ and its image by the bijection $\sigma$.}\label{fig:binary-tree+sigma}
\end{center}
\end{figure}

\titre{Tamari lattice.}
The Tamari lattice has a simple interpretation in terms of binary trees. The set of binary trees can be defined recursively by the following grammar. A binary tree $B$ is either a leaf denoted by~$\circ$ or is an ordered pair of binary trees, denoted $B=(B_1,B_2)$. 
It is often convenient to draw a binary tree by representing the leaf by a white vertex and the tree $B=(B_1,B_2)$ by a black vertex at the bottom joined to the subtrees $B_1$ (on the left) and $B_2$ (on the right). The tree $(((\circ,\circ),((\circ,\circ),\circ)),(\circ,(\circ,\circ)))$ is represented in Figure~\ref{fig:binary-tree+sigma}.

The set $\BB_n$ of binary trees with $n$ nodes has cardinality $C_n=\frac{1}{n+1}{2n \choose n}$ and there are well known bijections between the set $\BB_n$ and the set $\DD_n$. We call $\sigma$ the bijection defined as follows: the image of the binary tree reduced to a leaf is the empty word and the image of the binary tree $B=(B_1,B_2)$  is the Dyck path $\sigma(B)=\sigma(B_1)N\sigma(B_2)S$. An example is given in Figure \ref{fig:binary-tree+sigma}. \\

\begin{figure}[htb!]
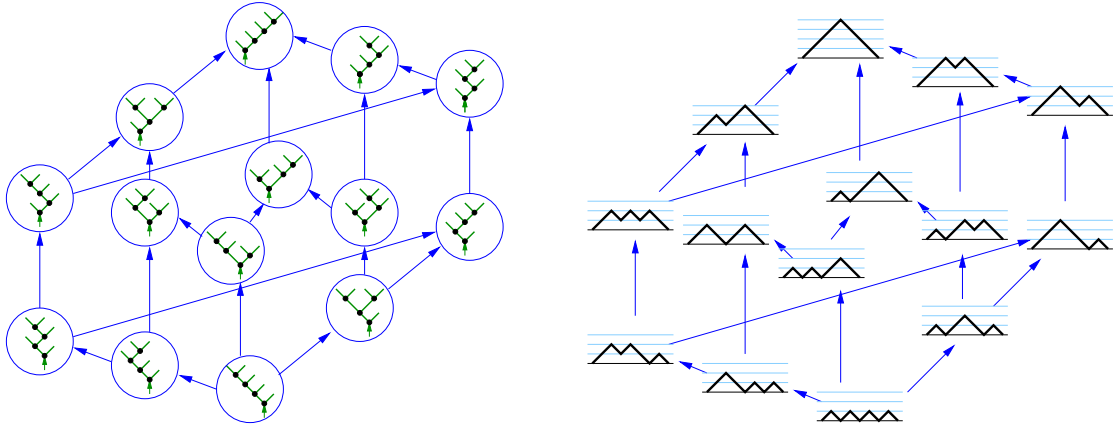

\begin{center}
\input{Tamari-lattice4.pstex_t}\hspace{1cm}\input{Tamari-lattice4-Dyck2.pstex_t}
\caption{Hasse diagram of the Tamari lattice $\LT_4$.}\label{fig:Tamari-lattice}
\end{center}
\end{figure}

In \cite{Tamari:lattice}, Tamari defined a partial order on the set $\BB_n$ of binary trees and proved to be a lattice. The covering relation for the Tamari lattice is defined has follows: a binary tree~$B$ containing a subtree of type $X=((B_1,B_2),B_3)$ is covered by the binary tree $B'$ obtained from $B$ by replacing $X$ by  $(B_1,(B_2,B_3))$. The Hasse diagram of the Tamari lattice on the set of binary trees with $4$ nodes is represented in Figure~\ref{fig:Tamari-lattice}~(left). 

The bijection $\sigma$ allows to transfer the Tamari lattice to the set of $\DD_n$ Dyck paths. We denote by $\LT_n$ the image of the Tamari lattice on  $\DD_n$ and denote by $P\lT Q$ if the path $P$ is less than or equal to the path $Q$ for this order. The  Hasse diagram of $\LT_4$ is represented in Figure~\ref{fig:Tamari-lattice}~(right).  The following ptoposition expresses the covering relation of the Tamari lattice  $\LT_n$ in terms of Dyck paths. This covering relation is illustrated in Figure~\ref{fig:covering}~$(b)$.

\begin{prop}\label{lem:covering-Tamari}
Let $P=NS^{\al_1}\ldots NS^{\al_n}$ and $Q=NS^{\be_1}\ldots NS^{\be_n}$ be two Dyck paths. The path $P$ is covered by the path $Q$ in the Tamari lattice $\LT_n$ if $Q$ is obtained from $P$ by swapping an $S$ step and the prime Dyck subpath following it, that is, there are indices $1\leq i<j\leq n$ with $\al_i>0$ and $i\underder{P} j$ such that $\be_i=\al_i-1$, $\be_j=\al_j+1$ and $\be_k=\al_k$ for all $k\neq i,j$.
\end{prop}

\begin{cor}\label{lem:Stanley-Tamari}
The Stanley lattice $\LS_n$ is a refinement of the Tamari lattice $\LT_n$. That is, for any pair of Dyck paths $P,Q$, $P\lT Q$ implies $P\lS Q$.
\end{cor}

\demof{Proposition \ref{lem:covering-Tamari}} Let $B$ be a binary tree and let $P=\sigma(B)$. \\
\ite We use the well known fact that \emph{there is a one-to-one correspondence between the subtrees of $B$ and the Dyck subpaths of $P$ which are either a prefix of $P$ or are preceded by an $N$ step}. (This classical property is easily shown by induction on the size of $P$.)\\
\ite If  the binary tree $B'$ is obtained from $B$ by replacing a subtree $X=((B_1,B_2),B_3)$ by $X'=(B_1,(B_2,B_3))$, then  the Dyck path $Q=\sigma(B')$ is obtained from $P$ by replacing a subpath $\sigma(X)=\sigma(B_1)N\sigma(B_2)SN\sigma(B_3)S$ by $\sigma(X')=\sigma(B_1)N\sigma(B_2)N\sigma(B_3)SS$; hence by swapping an $S$ step and the prime Dyck subpath following it.\\
\ite Suppose conversely that the Dyck path $Q$ is obtained from $P$ by swapping an $S$ step with a prime Dyck subpath $NP_3S$ following it. Then, there are two Dyck paths $P_1$ and $P_2$ (possibly empty) such that $W=P_1NP_2SNP_3S$ is a Dyck subpath of $P$ which is either a prefix of $P$ or is preceded by an $N$ step. Hence, the binary tree $B$ contains the subtree  $X=\sigma^{-1}(W)=((B_1,B_2),B_3)$, where $B_i=\sigma^{-1}(P_i),~i=1,2,3$. Moreover, the binary tree $B'=\sigma^{-1}(Q)$ is obtained from $B$ by replacing the subtree  $X=((B_1,B_2),B_3)$ by  $X'=(B_1,(B_2,B_3))=\sigma^{-1}(P_1NP_2NP_3SS)$.
\findem

\begin{figure}[ht]
\begin{center}
\input{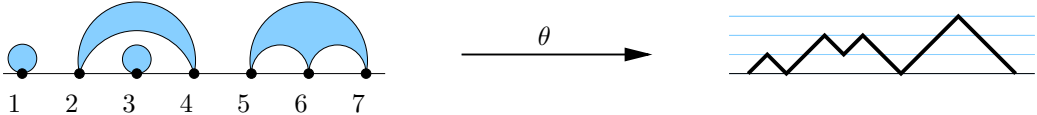}
\caption{A non-crossing partition and its image by the bijection $\theta$.}\label{fig:non-crossing+theta}
\end{center}
\end{figure}

\titre{Kreweras lattice.}
A partition of $\{1,\ldots,n \}$ is \emph{non-crossing} if whenever four elements $1\leq i<j<k<l\leq n$ are such that $i,k$ are in the same class and $j,l$ are in the same class, then the two classes coincide. The non-crossing partition whose classes are $\{1\}$, $\{2,4\}$,  $\{3\}$, and  $\{5,6,7\}$ is represented in Figure~\ref{fig:non-crossing+theta}. In this figure, each class is represented by a connected cell incident to the integers it contains.\\

The set $\NC_n$ of non-crossing partition on $\{1,\ldots,n \}$ has cardinality $C_n=\frac{1}{n+1}{2n \choose n}$ and there are well known bijections between non-crossing partitions and Dyck paths. We consider the bijection $\theta$ defined as follows. The image of a non-crossing partition $\pi$ of size $n$ by the mapping $\theta$ is the Dyck path $\theta(\pi)=NS^{\al_1}NS^{\al_2}\ldots NS^{\al_n}$, where $\al_i$ is the size of the class containing $i$ if $i$ is maximal in its class and $\al_i=0$ otherwise. An example is given in Figure  \ref{fig:non-crossing+theta}.\\

\begin{figure}[htb!]
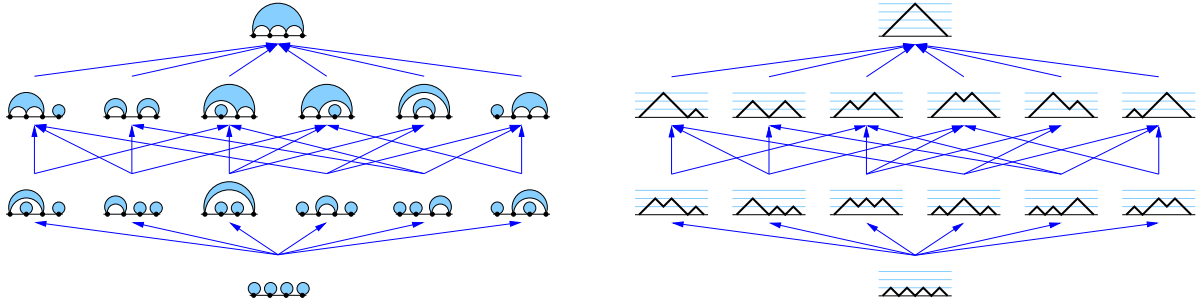

\begin{center}
\input{Kreweras-lattice4.pstex_t}\hspace{1cm}\input{Kreweras-lattice4-Dyck.pstex_t}
\caption{Hasse diagram of the Kreweras lattice $\LK_4$.}\label{fig:Kreweras-lattice}
\end{center}
\end{figure}

In \cite{Kreweras:non-crossing}, Kreweras showed that the partial order of refinement defines a lattice on the set $\NC_n$ of non-crossing partitions. The covering relation of this lattice corresponds to the merging of two parts when this operation does not break the \emph{non-crossing condition}. The Hasse diagram of the Kreweras lattice on the set $\NC_4$  is represented in Figure~\ref{fig:Kreweras-lattice}~(left). \\

\begin{figure}[ht]
\begin{center}
\input{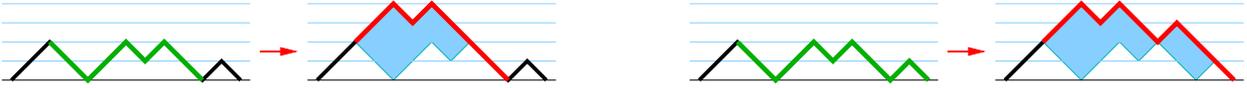}
\caption{Two examples of covering relation in the Kreweras lattice.}\label{fig:covering-Kreweras}
\end{center}
\end{figure}

The bijection $\theta$ allows to transfer the Kreweras lattice on the set $\DD_n$ of Dyck paths. We denote by $\LK_n$ the lattice structure obtained on $\DD_n$ and denote by $P\lK Q$ if the path $P$ is less than or equal to the path $Q$ for this order. The  Hasse diagram of $\LK_4$ is represented in Figure~\ref{fig:Kreweras-lattice}~(right). The following proposition expresses the covering relation of the Kreweras lattice  $\LK_n$ in terms of Dyck paths. This covering relation is represented in Figure~\ref{fig:covering-Kreweras}.

\begin{prop}\label{lem:covering-Kreweras} 
Let $P=NS^{\al_1}\ldots NS^{\al_n}$ and $Q=NS^{\be_1}\ldots NS^{\be_n}$ be two Dyck paths of size $n$. The path $P$ is covered by the path $Q$ in the Kreweras lattice $\LK_n$ if  $Q$ is obtained from $P$ by swapping a (non-empty) descent with a Dyck subpath following it, that is,  there are indices $1\leq i<j\leq n$ with $\al_i>0$ and $i\under{P} j$ such that $\be_i=0$, $\be_j=\al_i+\al_j$ and $\be_k=\al_k$ for all $k\neq i,j$.
\end{prop}

\begin{cor}\label{cor:Tamari-Kreweras}
The Tamari lattice $\LT_n$ is a refinement of the Kreweras lattice $\LK_n$. That is, for any pair $P,Q$ of Dyck paths, $P\lK Q$ implies $P\lT Q$.
\end{cor}

Proposition \ref{lem:covering-Kreweras} is a immediate consequence of the following lemma.

\begin{lemma}
Let $\pi$ be a non-crossing partition and let $P=\theta(\pi)$. Let $c$ and $c'$ be two classes of $\pi$ with the convention that  $i=\max(c)< j=\max(c')$. Then, the classes $c$ and $c'$ can be merged without breaking the non-crossing condition if and only if $i \under{P} j$. 
\end{lemma}

\dem For any index $k=1,\ldots,n$, we denote by $c_k$ the class of $\pi$ containing $k$. Observe that the classes $c$ and $c'$ can be merged without breaking the non-crossing condition if and only if there are no integers $r,s$ with $c_r=c_s$ such that  $r<i<s<j$ or  $i<r<j<s$. Observe also from the definition of the mapping~$\theta$ that for all index $l=1,\ldots,n$, the exceedence $e_l(P)$ is equal to the number of indices $k\leq l$ such that $\max(c_k)>l$.  \\
\ite We suppose that $i\under{P} j$ and we want to prove that merging the classes $c$ and $c'$ does not break the non-crossing condition. We first prove that there are no integers $r,s$ such that $i<r<j<s$ and $c_r= c_s$. Suppose the contrary. In this case,  there is no integer $k\leq r$ such that $r<\max(c_k)\leq j$ (otherwise, $c_k=c_r=c_s$ by the non-crossing condition, hence $\max(c_k)\geq \max(c_s)>j$).  Thus, $\{k\leq r/\max(c_k)>r\}= \{k\leq r/\max(c_k)>j\}\subsetneq  \{k\leq j/\max(c_k)>j\}$. This implies $e_r(P)<e_j(P)$ and contradicts the assumption $i\under{P}j$.  It remains to prove that there are no integers $r,s$ such that $r<i<s<j$ and $c_r= c_s$. Suppose the contrary and let $s'=\max(c_r)$. The case where $s'\geq j$ has been treated in the preceding point so we can assume that $s'<j$. In this case, there is no integer $k$ such that $i<k\leq s'$ and $\max(c_k)>s'$ (otherwise, $c_k=c_r=c_{s'}$ by the non-crossing condition, hence $\max(c_k)=\max(c_r)=s'$).  Thus, $\{k\leq i/\max(c_k)>i\} \subsetneq \{k\leq i/\max(c_k)>s'\}=\{k\leq s'/\max(c_k)>s'\}$. This implies $e_i(P)<e_{s'}(P)$ and contradicts the assumption $i\under{P}j$. \\
\ite We suppose now that merging the classes $c$ and $c'$ does not break the non-crossing partition and we want to prove that $i\under{P} j$. Observe that there is no integer $k$ such that $i<k\leq j$ and $\max(c_k)>j$ (otherwise, merging the classes $c$ and $c'$ would break the non-crossing condition). Thus, $\{k\leq j/\max(c_k)>j\} = \{k\leq i/\max(c_k)>j\}  \subseteq \{k\leq i/\max(c_k)>i\}$. This implies $e_j(P)\leq e_i(P)$. It remains to prove that there is no index $s$ such that $i<s<j$ and $e_s(P)<e_i(P)$. Suppose the contrary and consider the minimal such $s$. Observe that $s$ is maximal in its class, otherwise $e_{s-1}(P)=e_s(P)-1<e_i(P)$ contradicts the minimality of $s$. Observe also that $i<r=\min(c_s)$ otherwise merging the classes $c$ and $c'$ would break the non-crossing condition. By the non-crossing condition, there is no integer $k<r$ such that $r\leq \max(c_k)\leq s$. Thus,  $ \{k\leq r-1/\max(c_k)>r-1\} = \{k\leq r-1/\max(c_k)>s\}  \subseteq \{k\leq s/\max(c_k)>s\}$. This implies $e_{r-1}(P)\leq e_{s}(P)<e_i(P)$ and contradicts the minimality of $s$.\findem

\section{A bijection between Stanley intervals and realizers}\label{section:bijection}
In this section, we recall some definitions about triangulations and realizers. Then, we define a bijection between pairs of non-crossing Dyck paths and realizers.\\

\subsection{Triangulations and realizers}~\\
\textbf{Maps.} A \emph{planar map}, or \emph{map} for short, is an embedding of a connected finite planar graph in the sphere considered up to continuous deformation.  In this paper, maps have no loop nor multiple edge. The \emph{faces} are the connected components of the complement of the graph. By removing the midpoint of an edge we get two \emph{half-edges}, that is, one dimensional cells incident to one vertex. Two consecutive half-edges around a vertex define a \emph{corner}. If an edge is oriented we call \emph{tail} (resp. \emph{head}) the half-edge incident to the origin (resp. end). \\

\begin{figure}[htb!]
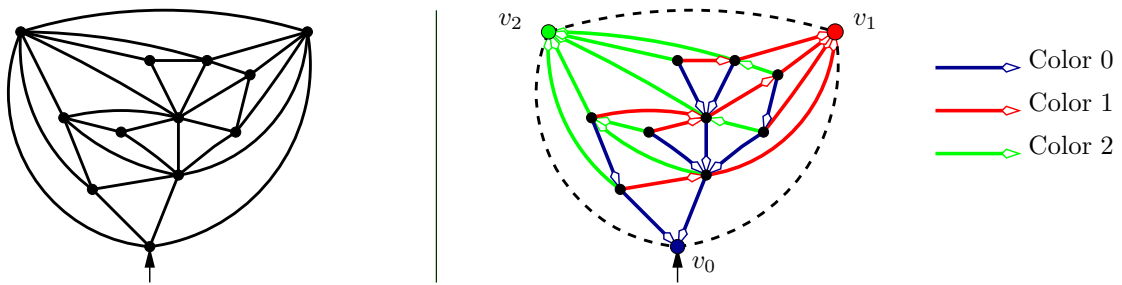

\begin{center}
\input{rooted-triangulation.pstex_t}\hspace{.7cm}\input{bijection-final.pstex_t}
\caption{A rooted triangulation (left) and one of its realizers (right).}\label{fig:rooted-triangulation}
\end{center}
\end{figure}

A \emph{rooted} map is a map together with a special half-edge which is not part of a complete edge and is called the \emph{root}. (Equivalently, a rooting is defined by the choice of a corner.) The root is incident to one vertex called \emph{root-vertex} and one face (containing it) called the \emph{root-face}.  When drawing maps in the plane the root is represented by an arrow pointing on the root-vertex and the root-face is the infinite one. See Figure~\ref{fig:rooted-triangulation} for an example. The vertices and edges incident to the root-face are called \emph{external} while the others are called \emph{internal}. From now on, \emph{maps are rooted} without further notice.\\

\titre{Triangulations.} A \emph{triangulation} is a map in which any face \emph{has degree} 3 (has 3 corners). A triangulation has \emph{size} $n$ if it has  $n$ internal vertices. The incidence relation between faces and edges together with Euler formula show that  a triangulation of size $n$ has $3n$ internal edges and $2n+1$ internal triangles. \\


In one of its famous \emph{census} paper, Tutte proved by a generating function approach that the number of triangulations of size $n$ is $t_n=\frac{2(4n+1)!}{(n+1)!(3n+2)!}$ \cite{Tutte:census1}.
A bijective proof of this result was given in \cite{Schaeffer:triang-3-connected}.\\

\titre{Realizers.} We now recall the notion of \emph{realizer} (or \emph{Schnyder wood}) defined by Schnyder \cite{Schnyder:wood1,Schnyder:wood2}. Given an edge coloring of a map, we shall call \emph{$i$-edge} (resp. \emph{$i$-tail}, \emph{$i$-head}) an edge (resp. tail, head) of color~$i$.

\begin{figure}[ht]
\begin{center}
\input{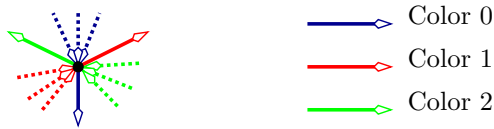}
\caption{Edges coloration and orientation around a vertex in a realizer (Schnyder condition).}
\label{fig:neighbors}
\end{center}
\end{figure}

\begin{Def}[\cite{Schnyder:wood1}] \label{def:realizer}
Let $M$ be a triangulation and let $U$ be the set of its internal vertices. Let $v_0$ be the root-vertex and let $v_1$, $v_2$ be the other external vertices with the convention that $v_0$, $v_1$, $v_2$ appear in counterclockwise order around the root-face.\\
A \emph{realizer} of $M$ is a coloring of the internal edges in three colors $\{0,1,2\}$ such that:
\begin{enumerate}
\item \emph{Tree condition}: for $i=0,1,2$, the $i$-edges form a tree $T_i$ with vertex set $U\cup \{v_i\}$. The vertex $v_i$ is considered to be the root-vertex of $T_i$ and the $i$-edges are oriented toward $v_i$.
\item \emph{Schnyder condition}: in clockwise order around any internal vertex there is: one 0-tail, some 1-heads, one 2-tail, some 0-heads, one 1-tail, some 2-heads. This situation is represented in Figure~\ref{fig:neighbors}. 
\end{enumerate}
We denote by $R=(T_0,T_1,T_2)$ this realizer. 
\end{Def}

A realizer is represented in Figure~\ref{fig:rooted-triangulation}~(right). Let $R=(T_0,T_1,T_2)$ be a realizer. We denote by $\bbar{T_0}$ the tree made of $T_0$ together with the edge $(v_0,v_1)$. For any internal vertex $u$, we denote by $\Pa_i(u)$ the parent of $u$ in the tree $T_i$. 
A \emph{cw-triangle} (resp. \emph{ccw-triangle}) is a triple of vertices $(u,v,w)$ such that $\Pa_0(u) = v, \Pa_2(v)=w$ and $\Pa_1(w) = u$ (resp. $\Pa_0(u) = v, \Pa_1(v)=w$ and $\Pa_2(w) = u$). A realizer is  called \emph{minimal} (resp. \emph{maximal}) if it has no cw-triangle (resp. ccw-triangle). It was proved in \cite{Mendez:these,Propp:lattice} that every triangulation has a unique minimal (resp. maximal) realizer. (The appellations \emph{minimal} and \emph{maximal} refer to a classical lattice which is defined on the set of realizers of any given triangulation \cite{Mendez:these,Propp:lattice}.)\\




\subsection{A bijection between pairs of non-crossing Dyck paths and realizers}~\\
In this subsection, we give an alternative (and simpler) description of the bijection defined  in \cite{Bonichon:realizers} between realizers and pairs of non-crossing Dyck paths.\\

\begin{figure}[h!]
\begin{center}
\input{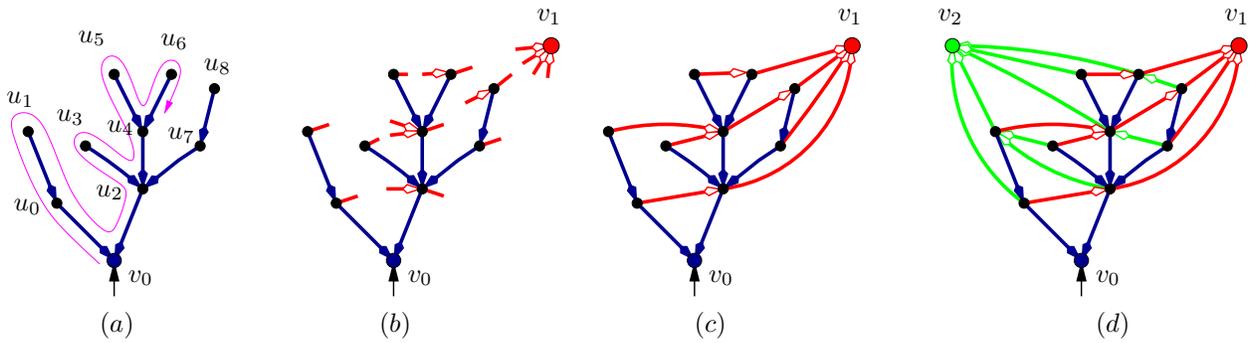}
\caption{Steps of the mapping $\Phi:(P,Q)\mapsto (T_0,T_1,T_2)$. (a) Step 1: build the tree $T_0$. (b) Step 2: add the  1-tails and 1-heads. (c) Step 3: join the 1-tails and 1-heads together. (d) Step 4: determine the third tree $T_2$.}
\label{fig:bijection}
\end{center}
\end{figure}

We first recall a classical bijection between plane trees and Dyck paths. A \emph{plane tree} is a rooted map whose underlying graph is a tree. Let $T$ be a plane tree. We \emph{make the tour} of the tree $T$ by following its border in clockwise direction starting and ending at the root (see Figure~\ref{fig:bijection}~(a)). We denote by $\omega(T)$ the word obtained by making the tour of the tree $T$ and writing $N$ the first time we follow an edge and $S$ the second time we follow this edge.  For instance, $w(T)=NNSSNNSNNSNSSNNSSS$ for the tree in Figure~\ref{fig:bijection}~(a). It is well known that the mapping $\omega$ is a bijection between plane trees with $n$ edges and Dyck paths of size~$n$ \cite{Knuth:generating-all-trees}.\\


Let $T$ be a plane tree. Consider the order in which the vertices are encountered while making the tour of $T$. This defines the \emph{clockwise order around} $T$ (or \emph{preorder}). For the tree in Figure~\ref{fig:bijection}~(a), the clockwise order is $v_0< u_0< u_1 <\ldots < u_8$. The tour of the tree also defines an order on the set of corners around each vertex $v$. We shall talk about the \emph{first} (resp. \emph{last}) \emph{corner of $v$ around $T$}.\\

We are now ready to define a mapping $\Psi$ which associates an ordered pair of Dyck paths to each realizer. 
\begin{Def}\label{def:Phi}
Let $M$ be a rooted triangulation of size $n$ and let $R=(T_0,T_1,T_2)$ be a realizer of $M$. Let $u_0,u_1,\ldots,u_{n-1}$ be the internal vertices of $M$ in clockwise order around $T_0$. Let $\be_i, i=1,\ldots,n-1$ be the number of 1-heads incident to $u_i$ and let $\be_n$ be the number of 1-heads incident to $v_1$. Then $\Psi(R)=(P,Q)$, where $P=\omega^{-1}(T_0)$ and $Q=NS^{\be_1}\ldots NS^{\be_n}$. 
\end{Def}

\begin{figure}[h!]
\begin{center}
\input{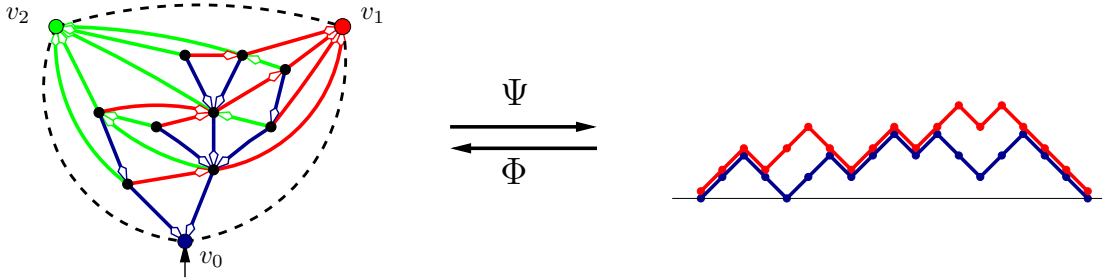}
\caption{The bijections $\Psi$ and $\Phi$.}\label{fig:Phi-and-Psi}
\end{center}
\end{figure}

The image of a realizer by the mapping $\Psi$ is represented in Figure~\ref{fig:Phi-and-Psi}.\\

\begin{thm}\label{thm:stanley}
The mapping $\Psi$ is a bijection between realizers of size $n$ and pairs of non-crossing Dyck paths of size $n$.
\end{thm}

The rest of this section is devoted to the proof of Theorem \ref{thm:stanley}. We first prove that the image of a realizer is indeed a pair of non-crossing Dyck paths.
\begin{prop}\label{prop:image-psi}
Let $R=(T_0,T_1,T_2)$ be a realizer of size $n$ and let $(P,Q)=\Psi(R)$.  Then, $P$ and $Q$ are both Dyck paths and moreover the path $P$ stays below the path $Q$. 
\end{prop}

Proposition \ref{prop:image-psi} is closely related to the Lemma \ref{lem:T1-left-to-right} below which, in turn, relies on the following technical lemma.
\begin{lemma}\label{lem:pigeon-hole}
Let $M$ be a map in which every face has degree three. We consider an orientation of the internal edges of $M$ such that every internal vertex has outdegree 3 (i.e. is incident to exactly 3 tails). Let $C$ be a simple cycle made of $c$ edges. By the Jordan Lemma, the cycle $C$ separates the sphere into two connected regions. We call \emph{inside} the region not containing the root. Then, the number of tails incident with $C$ and lying strictly inside $C$ is $c-3$.
\end{lemma}

\dem 
Let $v$ (resp. $f,~e$) be the number of vertices (resp. faces, edges) lying strictly inside $C$. Note that the edges strictly inside $C$ are internal hence are oriented. The number $i$ of tails incident with $C$ and lying strictly inside $C$ satisfies $e=3v+i$. Moreover, the incidence relation between edges and faces implies $3f=2e+c$ and the Euler relation implies $(f+1)+(v+c)=(e+m)+2$. Solving for $i$ gives $i=c-3$.
\findem

\begin{lemma}\label{lem:T1-left-to-right}
Let $R=(T_0,T_1,T_2)$ be a realizer. Then, for any 1-edge $e$ the tail of $e$ is encountered before its head around the tree $\bbar{T_0}$.
\end{lemma}

\begin{figure}[h!]
\begin{center}
\input{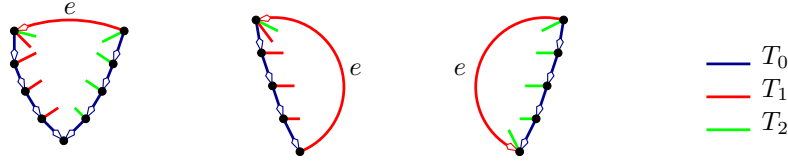}
\vspace{-.3cm}
\caption{Case analysis for a 1-edge $e$ whose head appears before its tail around the tree $\bbar{T_0}$.}\label{fig:reverse-1-edge}
\end{center}
\end{figure}

\demof{Lemma  \ref{lem:T1-left-to-right}}
Suppose a 1-edge $e$ breaks this rule and consider the cycle $C$ made of $e$ and the 0-path joining its endpoints. Using the Schnyder condition it is easy to show that the number of tails incident with $C$ and lying strictly inside $C$ is equal to the number of edges of $C$ (the different possibilities are represented in Figure \ref{fig:reverse-1-edge}). This contradicts Lemma \ref{lem:pigeon-hole}.
\findem

\vspace{-.3cm}

\begin{lemma}\label{lem:tail-head-sequence}
Let $P=NS^{\al_1}\ldots NS^{\al_n}$ be a Dyck path and let $T=\omega^{-1}(P)$. Let $v_0$ be the root-vertex of the tree $T$ and  let $u_0,u_1,\ldots,u_{n-1}$ be its other vertices in clockwise order around $T$. Then, the word obtained by making the tour of $T$ and writing  $\BS^{\be_i}$  when arriving at the first corner of  $u_i$ and $\BN$ when arriving at the last corner of $u_i$ is $\textbf{W}=\BS^{\be_0}\BN^{\al_1}\BS^{\be_1}\ldots \BS^{\be_{n-1}}\BN^{\al_n}$.
\end{lemma}

\dem
We consider the word $\mathcal{W}$ obtained by making the tour of $T$ and writing $N\BS^{\be_i}$  when arriving at the first corner of  $u_i$ and $\BN S$ when arriving at the last corner of $u_i$ for $i=0,\ldots,n-1$. 
By definition of the mapping $\omega$, the restriction of $\mathcal{W}$ to the letters $N,S$ is $W=\omega(T)=P= NS^{\al_1}\ldots NS^{\al_n}$. Therefore, $\mathcal{W}=N\BS^{\be_0}(\BN S)^{\al_1}N\BS^{\be_1}(\BN S)^{\al_2} \ldots N\BS^{\be_{n-1}}(\BN S)^{\al_n}$. Hence, the restriction of $\mathcal{W}$ to the letters $\BN,~\BS$ is $\textbf{W}=\BS^{\be_0}\BN^{\al_1}\BS^{\be_1}\BN^{\al_2} \ldots \BS^{\be_{n-1}}\BN^{\al_n}$. 
\hfill \findembis{0cm}

\titre{Proof of Proposition~\ref{prop:image-psi}:} We denote $P=NS^{\al_1}\ldots NS^{\al_n}$ and $Q=NS^{\be_1}\ldots NS^{\be_n}$.\\
\ite The mapping $\omega$ is known to be a bijection between trees and Dyck paths, hence $P=\omega(T)$ is a Dyck path.\\
\ite We want to prove that $Q$ is a Dyck path staying above $P$. Consider the word $\textbf{W}$ obtained by making the tour of $\bbar{T_0}$ and writing $\BN$ (resp. $\BS$) when we encounter a 1-tail (resp. 1-head). By Lemma \ref{lem:tail-head-sequence}, the word $\textbf{W}$ is $\BS^{\be_0}\BN^{\al_1}\BS^{\be_1}\BN^{\al_2} \ldots \BS^{\be_{n-1}}\BN^{\al_n}\BS^{\be_n}$. By Lemma \ref{lem:T1-left-to-right}, the word $\textbf{W}$ is a Dyck path. In particular, $\BS^{\be_0}=0$ and $\sum_{i=1}^n \be_i =\sum_{i=1}^n \al_i=n$, hence the path $Q$ returns to the origin. Moreover, for all $i=1,\ldots,n$, $\delta_i(P,Q)=\sum_{j=1}^{n}\al_i-\be_i\geq 0$. Thus, the path $Q$ stays above $P$. In particular, $Q$ is a Dyck path.
\findem

In order to prove Theorem~\ref{thm:stanley}, we shall now define a mapping $\Phi$ from pairs of non-crossing Dyck paths to realizers and prove it to be the inverse of $\Psi$. We first define \emph{prerealizers}. 

\begin{Def} \label{def:prerealizer}
Let $M$ be a map. Let $v_0$ be the root-vertex, let $v_1$ be another external vertex and let $U$ be the set of the other vertices.
A \emph{prerealizer} of $M$ is a coloring of the edges in two colors $\{0,1\}$ such that:
\begin{enumerate}
\item \emph{Tree condition}: for $i=0,1$, the $i$-edges form a tree $T_i$ with vertex set $U\cup \{v_i\}$. The vertex $v_i$ is considered to be the root-vertex of $T_i$ and the $i$-edges are oriented toward $v_i$.
\item \emph{Corner condition}:  in clockwise order around any  vertex $u\in U$  there is: one 0-tail, some 1-heads, some 0-heads, one 1-tail.
\item \emph{Order condition}: for any 1-edge $e$, the tail of $e$ is encountered before its head around the tree $\bbar{T_0}$, where $\bbar{T_0}$ is the tree obtained from $T_0$ by adding the edge $(v_0,v_1)$ at the right of the root.
\end{enumerate}
We denote by $\PR=(T_0,T_1)$ this prerealizer. 
\end{Def}

An example of prerealizer is given in Figure~\ref{fig:bijection}~(c). 

\begin{lemma} \label{lem:trees2}
Let $\PR=(T_0,T_1)$ be a prerealizer. Then, there exists a unique tree $T_2$ such that $R=(T_0,T_1,T_2)$ is a realizer.
\end{lemma}

In order to prove Lemma  \ref{lem:trees2}, we need to study the sequences of corner around the faces of prerealizers.  If $h$ and $h'$ are two consecutive half-edges in clockwise order around a vertex $u$ we denote by $c=(h,h')$ the corner delimited by $h$ and $h'$. For $0\leq i,j\leq 2$, we call  \emph{$(h_i,h_j)$-corner} (resp. \emph{$(h_i,t_j)$-corner}, \emph{$(t_i,h_j)$-corner}, \emph{$(t_i,t_j)$-corner}) a corner $c=(h,h')$ where $h$ and $h'$ are respectively an $i$-head (resp. $i$-head, $i$-tail, $i$-tail) and a $j$-head (resp. $j$-tail, $j$-head, $j$-tail).\\

\demof{Lemma  \ref{lem:trees2}}
Let $\PR=(T_0,T_1)$ be a prerealizer and let $N=T_0\cup T_1$ be the underlying map. Let~$v_0$ (resp. $v_1$) be the root-vertex of $T_0$ (resp. $T_1$) and let $U$ be the set of vertices distinct from $v_0,v_1$. Let~$\bbar{T}_0$ (resp. $\bbar{N}$) be the tree (resp. map) obtained from $T_0$ (resp. $N$) by adding the edge $(v_0,v_1)$ at the right of the root. We first prove that there is at most one tree $T_2$ such that $R=(T_0,T_1,T_2)$ is a realizer. 
\begin{itemize}
\item Let $f$ be an internal face of $\bbar{N}$ and let $c_1,c_2,\ldots,c_k$ be the corners of $f$ encountered in clockwise order around $\bbar{T}_0$. Note that $c_1,c_2,\ldots,c_k$ also correspond to the clockwise order of the corners around the face $f$.
We want to prove the following properties: \\
\indent \iten the corner $c_1$ is a  $(t_1,t_0)$-corner,\\
\indent \iten the corner $c_2$ is either a $(h_0,h_0)$- or a $(h_0,t_1)$-corner,\\
\indent \iten the corners $c_3,\ldots,c_{k-1}$ are $(h_1,h_0)$-, $(h_1,t_1)$-, $(t_0,h_0)$- or $(t_0,t_1)$-corners,\\
\indent \iten the corner $c_k$ is either a $(h_1,h_1)$- or a $(t_0,h_1)$-corner.\\
First note that by the \emph{corner condition} of the prerealizers the possible corners are of type $(h_0,h_0)$, $(h_0,t_1)$, $(h_1,h_0)$, $(h_1,h_1)$, $(h_1,t_1)$, $(t_0,h_0)$, $(t_0,h_1)$, $(t_0,t_1)$ and $(t_1,t_0)$. By the \emph{order condition}, one enters a face for the first time  (during a tour of $T_0$) when crossing a 1-tail. Hence, the  first corner $c_1$ of $f$ is a $(t_1,t_0)$-corner while the corners $c_i,~i=2,\ldots,k$ are not $(t_1,t_0)$-corners.  Since $c_1$ is a  $(t_1,t_0)$-corner, the corner $c_2$ is either a  $(h_0,h_0)$- or a $(h_0,t_1)$-corner. Similarly,  since $c_1$ is a  $(t_1,t_0)$-corner, the corner $c_k$ is either a $(h_1,h_1)$- or a $(t_0,h_1)$-corner. Moreover, for $i=2,\ldots,k-1$, the corner $c_i$ is not a $(h_1,h_1)$- nor a $(t_0,h_1)$-corner or $c_{i+1}$ would be a $(t_1,t_0)$-corner. Therefore, it is easily seen by induction on $i$ that the corners $c_i,~i=3,\ldots,k-1$ are either $(h_1,h_0)$-,  $(h_1,t_1)$-, $(t_0,h_0)$- or $(t_0,t_1)$-corners.
\item By a similar argument we prove that the corners of the external face of $\bbar{N}$ are $(h_1,h_0)$-, $(h_1,t_1)$-, $(t_0,h_0)$- or $(t_0,t_1)$-corners except for the corner incident to $v_0$ which is a $(h_0,h_0)$-corner and the corner incident to $v_1$ which is a $(h_1,h_1)$-corner.
\item Let $v_2$ be an isolated vertex in the external face of $N$. If a tree $T_2$ with vertex set $U\cup \{v_2\}$ is such that $R=(T_0,T_1,T_2)$ is a realizer, then there is one 2-tail in each  $(h_1,h_0)$-, $(h_1,t_1)$-, $(t_0,h_0)$- or $(t_0,t_1)$-corner of $\bbar{N}$ while the 2-heads are only incident to the $(t_0,t_1)$-corners and to the vertex $v_2$. By the preceding points, there is exactly one $(t_1,t_0)$ corner in each internal face and none in the external face. Moreover there is at most one way of connecting the 2-tails and the 2-heads in each face of $\bbar{N}$. Thus, there is at most one tree $T_2$ such that $R=(T_0,T_1,T_2)$ is a realizer.
\end{itemize}
We now prove that there exists a tree $T_2$ such that  $R=(T_0,T_1,T_2)$ is a realizer. Consider the colored map $(T_0,T_1,T_2)$ obtained by \\
\iten  adding an isolated vertex $v_2$ in the external face of $\bbar{N}$. \\
\iten  adding a 2-tail in each  $(h_1,h_0)$-, $(h_1,t_1)$-, $(t_0,h_0)$- and $(t_0,t_1)$-corner of $\bbar{N}$.\\ 
\iten  joining each 2-tail in an internal face $f$ (resp. the external face) to the unique $(t_0,t_1)$-corner of $f$ (resp. to~$v_2$).\\
We denote by $M=T_0\cup T_1 \cup T_2 \cup \{(v_0,v_1),(v_0,v_2),(v_1,v_2)\}$ the underlying map. 
\begin{itemize}
\item We first prove that the map $M=T_0\cup T_1 \cup T_2 \cup \{(v_0,v_1),(v_0,v_2),(v_1,v_2)\}$ is a triangulation.
Let $f$ be an internal face. By a preceding point, $f$  has exactly one $(t_1,t_0)$ corner $c$  and the  $(h_1,h_0)$-, $(h_1,t_1)$-, $(t_0,h_0)$- or $(t_0,t_1)$-corners are precisely the ones that are not consecutive with $c$ around $f$. Thus, the internal faces of $N$ are triangulated (split into sub-faces of degree 3) by the 2-edges. Moreover, the only corners of the external face of $\bbar{N}$ which are not of type  $(h_1,h_0)$, $(h_1,t_1)$, $(t_0,h_0)$ or $(t_0,t_1)$ are the (unique) corner around  $v_0$ and  the (unique) corner around  $v_1$. Hence the external face of $\bbar{N}$ is triangulated by the 2-edges together with the edges $(v_0,v_2)$ and $(v_1,v_2)$. Thus, every face of $M$ has degree 3. It only remains to prove that $M$ has no multiple edge. Since the faces of $M$ are of degree 3 and every internal vertex has outdegree 3, the hypothesis of  Lemma \ref{lem:pigeon-hole} are satisfied. By this lemma, there can be no multiple edge (this would create a cycle of length 2 incident to -1 tails!). Thus, the map $M$ has no multiple edge and is a triangulation.
\item We now prove that the coloring $R=(T_0,T_1,T_2)$ is a realizer of $M$. By construction, $R$ satisfies de \emph{Schnyder-condition}. Hence it only remains to prove that $T_2$ is a tree. Suppose there is a cycle $C$ of 2-edges. Since every vertex in $C$ is incident to one 2-tail, the cycle $C$ is directed. Therefore, the Schnyder condition proves that there are $c=|C|$ tails incident with $C$ and lying strictly inside $C$. This contradicts Lemma \ref{lem:pigeon-hole}. Thus, $T_2$ has no cycle. Since $T_2$ has $|U|$ edges and $|U|+1$ vertices it is a tree.
\end{itemize}
\vspace{-.2cm}
\findem

We are now ready to define a mapping $\Phi$ from pairs of non-crossing Dyck paths to realizers. This mapping is illustrated by Figure \ref{fig:bijection}. Consider a pair of Dyck paths $P=NS^{\al_1}\ldots NS^{\al_n}$ and $Q=NS^{\be_1}\ldots NS^{\be_n}$ such that $P$ stays below $Q$. The image of $(P,Q)$ by the mapping $\Phi$ is the realizer $R=(T_0,T_1,T_2)$ obtained as follows.\\

\titre{Step 1.} The tree $T_0$ is $\omega^{-1}(P)$.  We denote by $v_0$ its root-vertex and by $u_0,\ldots,u_n$ the other vertices in clockwise order around $T_0$. We denote by $\bbar{T_0}$ the tree obtained from $T_0$ by adding a new vertex $v_1$ and an edge $(v_0,v_1)$ at the right of the root. \vspace{.1cm}\\
\titre{Step 2.} We glue a 1-tail in the last corner of each vertex $u_i,i=0,\ldots,n-1$ and we glue $\be_i$ 1-heads in the first corner of each vertex $u_i,i=1,\ldots,n-1$ (if $u_i$ is a leaf we glue the 1-heads before the 1-tail in clockwise order around $u_i$). We also glue $\be_n$ 1-heads in the (unique) corner of $v_1$. This operation is illustrated by Figure~\ref{fig:bijection}~(b)\vspace{.1cm}.\\
\titre{Step 3.} We consider the sequence of 1-tails and 1-heads around  $\bbar{T_0}$. Let $W$ be the word obtained by making the tour of $\bbar{T_0}$ and writing $N$ (resp. $S$) when we cross a 1-tail (resp. 1-head). By Lemma \ref{lem:tail-head-sequence}, $W=N^{\al_1}S^{\be_1}\ldots N^{\al_n}S^{\be_n}$. Since the path $P$ stays below the path $Q$, we have $\delta_i(P,Q)=\sum_{j\leq i}\al_j-\be_j\geq 0$ for all $i=1,\ldots,n$, hence $W$ is a Dyck path. Thus, there exists a unique way of joining each 1-tail to a 1-head that appears after it around the tree $\bbar{T_0}$ so that the 1-edges do not intersect (this statement is equivalent to the well-known fact that there is a unique way of matching parenthesis in a well parenthesized word); we denote by $T_1$ the set of 1-edges obtained in this way. This operation is illustrated in Figure~\ref{fig:bijection}~(c).\vspace{.1cm}\\
\titre{Step 4.} The set $T_1$ of 1-edges is a tree directed toward $v_1$; see Lemma \ref{lem:T1-is-tree} below. Hence, by construction, $\PR=(T_0,T_1)$ is a prerealizer. By Lemma \ref{lem:trees2}, there is a unique tree $T_2$ such that $R=(T_0,T_1,T_2)$ is a realizer and we define $\Phi(P,Q)=R$.\\

In order to prove that \emph{step 4} of the bijection $\Phi$  is well defined, we need the following lemma.
\begin{lemma}\label{lem:T1-is-tree}
The set $T_1$ of 1-edges obtained at step 3 in the definition of $\Phi$ is a tree directed toward $v_1$ and spanning the vertices in $U_1=\{u_0,\ldots,u_{n-1},v_1\}$.
\end{lemma}

\dem  
\ite Every vertex in $U_1$ is incident to an edge in $T_1$ since there is a 1-tail incident to each vertex $u_i,~i=1,\ldots,n-1$ and at least one 1-head incident to $v_1$ since $\be_n>0$. \\
\ite We now prove that the tree $T_1$ has no cycle. Since every vertex in $U_1$ is incident to at most one 1-tail, any 1-cycle is directed. Moreover, if $e$ is a 1-edge directed from $u_i$ to $u_j$ then $i<j$ since the last corner of $u_i$ appears before the first corner of $u_j$ around $T_0$. Therefore, there is no directed cycle. \\
\ite Since $T_1$ is a set of $n$ edges incident to $n+1$ vertices and having no cycle, it is a tree. Since the only sink is $v_1$, the tree $T_1$ is directed toward $v_1$ (make an induction on the size of the oriented tree $T_1$ by removing a leaf).\hfill \findembis{0cm}

The mapping $\Phi$ is well defined and the image of any pair of non-crossing Dyck paths is a realizer. Conversely, by Proposition~\ref{prop:image-psi}, the image of any realizer by $\Psi$ is a pair of non-crossing Dyck paths. It is clear from the definitions that $\Psi \circ \Phi$ (resp. $\Phi \circ \Psi$) is the identity mapping on pairs of non-crossing Dyck paths (resp. realizers). Thus, $\Phi$ and $\Psi$ are inverse bijections between realizers of size $n$ and pairs of non-crossing Dyck paths of size $n$. This concludes the proof of Theorem~\ref{thm:stanley}.\hfill\findembis{0cm}

\section{Intervals of the Tamari lattice}\label{section:Tamari}
In the previous section, we defined a bijection $\Phi$ between pairs of non-crossing Dyck paths and realizers. Recall that the pairs of non-crossing Dyck paths correspond to the intervals of the Stanley lattice. In this section, we study the restriction of the bijection $\Phi$ to the intervals of the Tamari lattice. 

\begin{thm}\label{thm:tamari}
The bijection $\Phi$ induces a bijection between the intervals of the Tamari lattice $\LT_n$ and minimal realizers of size $n$.
\end{thm}

Since every triangulation has a unique minimal realizer, Theorem \ref{thm:tamari} implies that the mapping $\Phi'$ which associates with a Tamari interval $(P,Q)$ the triangulation underlying $\Phi(P,Q)$ is a bijection. This gives a bijective explanation to the relation between the number of Tamari intervals enumerated in \cite{Chapoton:Tamari-intervals} and the number of triangulations enumerated in \cite{Tutte:census1,Schaeffer:triang-3-connected}.
 
\begin{cor}
The number of intervals in the Tamari lattice $\LT_n$ is equal to the number $\frac{2(4n+1)!}{(n+1)!(3n+2)!}$ of triangulations of size $n$.
\end{cor}

\begin{figure}[h!]
\begin{center}
\input{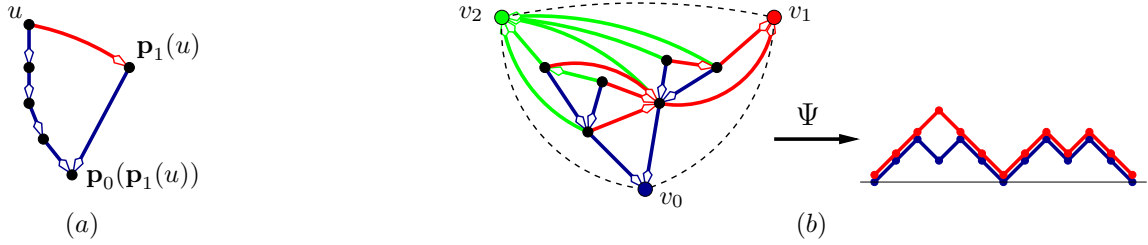}
\caption{(a) Characterization of minimality: $\Pa_0(\Pa_1(u))$ is an ancestor of $u$ in $T_0$. (b) A minimal realizer and its image by $\Psi$.}
\label{fig:minimality}
\end{center}
\end{figure}

The rest of this section is devoted to the proof of Theorem \ref{thm:tamari}. We first recall a characterization of minimality given in \cite{Bonichon:BGH02} and illustrated in Figure \ref{fig:minimality}.

\begin{prop}[\cite{Bonichon:BGH02}] \label{prop:branch}
A realizer $R=(T_0,T_1,T_2)$ is minimal if and only if for any internal vertex~$u$, the vertex $\Pa_0(\Pa_1(u))$ is an ancestor of $u$ in the tree $T_0$.
\end{prop}

Using Proposition \ref{prop:branch}, we obtain the following characterization of the pairs of non-crossing Dyck paths $(P,Q)$ whose image by the bijection $\Phi$ is a minimal realizer.

\begin{prop}\label{prop:characterization-min}
Let  $(P,Q)$ be a pair of non-crossing Dyck paths and let $R=(T_0,T_1,T_2)=\Phi(P,Q)$.  Let  $u_0,\ldots,u_{n-1}$ be the non-root vertices of $T_0$ in clockwise order. Then, the realizer $R$ is minimal if and only if $\delta_i(P,Q)\leq \delta_j(P,Q)$ whenever $u_i$ is the parent of $u_j$ in $T_0=\omega^{-1}(P)$.
\end{prop}


In order to prove Proposition \ref{prop:characterization-min}, we need to interpret the value of $\delta_i(P,Q)$ is terms of the realizer $R=\Phi(P,Q)$. Let $u$ be an internal vertex of the triangulation underlying the realizer  $R=(T_0,T_1,T_2)$. We say that a 1-tail is \emph{available at $u$} if this tail appears before the first corner of $u$ in clockwise order around $T_0$ while the corresponding 1-head appears (strictly) after the first corner of $u$.

\begin{lemma}\label{lem:nb-available-tails}
Let  $(P,Q)$ be a pair of non-crossing Dyck paths and let $R=(T_0,T_1,T_2)=\Phi(P,Q)$.  Let  $u_0,\ldots,u_{n-1}$ be the non-root vertices of $T_0$ in clockwise order. The number of 1-tails available at $u_i$ is~$\delta_i(P,Q)$.
\end{lemma}

\demof{Lemma \ref{lem:nb-available-tails}} We denote $P=NS^{\al_1}\ldots NS^{\al_n}$ and $Q=NS^{\be_1}\ldots NS^{\be_n}$.
Let $\mathcal{W}$ be the word obtained by making the tour of $T_0$ and writing $N\BS^{\be_i}$  when arriving at the first corner of  $u_i$ and $\BN S$ when arriving at the last corner of $u_i$ for $i=0,\ldots,n-1$ (with the convention that $\be_0=0$). By definition of the mapping $\omega$, the restriction of $\mathcal{W}$ to the letters $N,S$ is $\omega(T_0)=P= NS^{\al_1}\ldots NS^{\al_n}$. Therefore, $\mathcal{W}=N\BS^{\be_0}(\BN S)^{\al_1}N\BS^{\be_1}(\BN S)^{\al_2} \ldots N\BS^{\be_{n-1}}(\BN S)^{\al_n}$. The prefix of $\mathcal{W}$ written after arriving at the first corner of  $u_i$ is $N\BS^{\be_0}(\BN S)^{\al_1}N\BS^{\be_1}\ldots (\BN S)^{\al_i}N\BS^{\be_{i}}$. The sub-word $\BS^{\be_0}\BN^{\al_1}\BS^{\be_1}\ldots \BN^{\al_i}\BS^{\be_{i}}$ corresponds to the sequence of 1-tails and 1-heads encountered so far ($\BN$ for a 1-tail, $\BS$ for a 1-head).  Thus, the number of 1-tails available at  $u_i$ is $\sum_{j\leq i}\al_j-\be_j=\delta_i(P,Q)$.
\findem

\demof{Proposition \ref{prop:characterization-min}}\\
\ite We suppose that a vertex $u_i$ is the parent of a vertex $u_j$ in $T_0$ and that $\delta_i(P,Q)> \delta_j(P,Q)$, and we want to prove that the realizer $R=\Phi(P,Q)$ is not minimal.  Since $u_i$ is the parent of $u_j$ we have $i<j$ and all the vertices $u_r,~i<r\leq j$ are descendants of $u_i$. By Lemma \ref{lem:nb-available-tails},  $\delta_i(P,Q)> \delta_j(P,Q)$ implies that there is a 1-tail $t$ available at $u_i$  which is not available at $u_j$, hence the corresponding 1-head is incident to a vertex $u_l$ with $i<l\leq j$. Let $u_k$ be the vertex incident to the 1-tail~$t$. Since $t$ is available at $u_i$, the vertex $u_k$ is not a descendant of $u_i$. But $\Pa_0(\Pa_1(u_k))=\Pa_0(u_l)$ is either $u_i$ or a descendant of $u_i$ in $T_0$. Thus, the vertex $u_k$ contradicts the minimality condition given by Proposition \ref{prop:branch}. Hence, the realizer $R$ is not minimal.\\
\ite We suppose that  the realizer  $R$ is not minimal and we want to prove that there exists a vertex $u_i$ parent of a vertex $u_j$ in $T_0$ such that  $\delta_i(P,Q)>\delta_j(P,Q)$. By Proposition \ref{prop:branch}, there exists a vertex $u$ such that $\Pa_0(\Pa_1(u))$ is not an ancestor of $u$ in $T_0$. In this case, the 1-tail $t$ incident to $u$ is available at $u_i=\Pa_0(\Pa_1(u))$ but not at $u_j=\Pa_1(u)$ (since $t$ cannot appear between the  first corner of $u_i$  and the first corner of $u_j$ around $T_0$, otherwise $u$ would be a descendant of $u_i$). Moreover, any 1-tail $t'$ available at $u_j$ appears before the 1-tail $t$ around $T_0$ (otherwise, the 1-edge corresponding to $t'$ would cross the 1-edge $(u,u_j)$). Hence, any 1-tail $t'$ available at $u_j$ is also available at~$u_i$. Thus, there are more 1-tails available at $u_i$ than at $u_j$. By Lemma~\ref{lem:nb-available-tails}, this implies $\delta_i(P,Q)> \delta_j(P,Q)$.
\findem

\begin{prop}\label{prop:characterization-tamari}.
Let  $(P,Q)$ be a pair of non-crossing Dyck paths. Let $T=\omega^{-1}(P)$, let $v_0$ be the root-vertex of the tree $T$ and let $u_0,\ldots,u_{n-1}$ be its other vertices in clockwise order. Then, $P\lT Q$ if and only if  $\delta_i(P,Q)\leq \delta_j(P,Q)$ whenever $u_i$ is the parent of $u_j$.
\end{prop}

Propositions \ref{prop:characterization-min} and Propositions \ref{prop:characterization-tamari} clearly imply Theorem~\ref{thm:tamari}. Hence, it only remains to prove  Proposition \ref{prop:characterization-tamari}. \\

\dem We denote  $Q=NS^{\be_1}\ldots NS^{\be_n}$. \\
\ite We suppose that $P\lT Q$ and want to prove that $\delta_k(P,Q)\leq \delta_l(P,Q)$ whenever $u_k$ is the parent of $u_l$. We make an induction on $\Delta(P,Q)$. If $\Delta(P,Q)=0$, then $P=Q$ and the property holds. If $\Delta(P,Q)>0$ there is a path $Q'=NS^{\be_1'}\ldots NS^{\be_n'}$ such that $P\lT Q'$ and $Q'$ is covered by $Q$ in the Tamari lattice. The three paths $P,Q',Q$ are represented in Figure \ref{fig:characterize-Tamari}. By definition, there are two indices $1\leq i<j\leq n$ such that $i\underder{Q'}j$ and $\be_i=\be_i'+1$, $\be_j=\be_j-1$ and $\be_k=\be_k'$ for all $k\neq i,j$. Thus, $\delta_k(P,Q)=\delta_k(P,Q')+1$ if $i\leq k<j$ and  $\delta_k(P,Q)=\delta_k(P,Q')$ otherwise. By the induction hypothesis we can assume that  $\delta_k(P,Q')\leq \delta_l(P,Q')$ whenever $u_k$ is the parent of $u_l$. Suppose there exists $u_k$ parent of $u_l$ such that $\delta_k(P,Q)>\delta_l(P,Q)$. Note that if $u_k$ is the parent of $u_l$ then $k<l$ and for all $k<r\leq l$, the vertex $u_r$ is a proper descendant of $u_k$. Since $\delta_k(P,Q)>\delta_l(P,Q)$ and $\delta_k(P,Q')\leq \delta_l(P,Q')$ we have $k<j\leq l$, hence $u_j$ is a proper descendant of $u_k$. Note that for all $r=0,\ldots,n-1$,  $e_r(P)+1$ is equal to the height of the vertex $u_r$ in the tree $T$ (i.e. the distance between $v_0$ and $u_r$). Thus, $e_k(P)<e_j(P)$. Moreover, by the induction hypothesis, $\delta_k(P,Q')\leq \delta_j(P,Q')$. Hence, $e_k(Q')=e_k(P)+\delta_k(P,Q')<e_j(Q')=e_j(P)+\delta_j(P,Q')$. But since $i\leq k<j$ this contradicts the hypothesis  $i\underder{Q'}j$. We reach a contradiction, hence  $\delta_k(P,Q)\leq \delta_l(P,Q)$ whenever $u_k$ is the parent of $u_l$.
\begin{figure}[h!]
\begin{center}
\input{characterize-Tamari.pstex_t}
\caption{The Dyck paths $P\lT Q' \lT Q$.}\label{fig:characterize-Tamari}
\end{center}
\end{figure}

\ite We suppose that $\delta_k(P,Q)\leq \delta_l(P,Q)$ whenever $u_k$ is the parent of $u_l$ and want to prove that $P\lT Q$. We make an induction on $\Delta(P,Q)$. If $\Delta(P,Q)=0$, then $P=Q$ and the property holds. Suppose $\Delta(P,Q)>0$ and let $\delta=\max\{\delta_k(P,Q),k=0\ldots n\}$, let $e=\min\{e_k(P)/\delta_k(P,Q)=\delta\}$ and let $i=\max\{k/e_k(P)=e \textrm{ and } \delta_k(P,Q)=\delta\}$. Let $j$ be the first index such that $i<j\leq n$ and $u_j$ is not a descendant of $u_i$ ($j=n$ if $u_{i+1},\ldots,u_{n-1}$ are all descendants of $u_i$). Let $Q'=NS^{\be_1'}\ldots NS^{\be_n'}$ with  $\be_i'=\be_i+1$, $\be_j'=\be_j-1$ and $\be_k'=\be_k$ for all $k\neq i,j$. The paths $P,~Q$ and $Q'$ are represented in Figure \ref{fig:characterize-Tamari}. We want to prove that \emph{$Q'$ is a Dyck path covered by $Q$ in the Tamari lattice and $P\lT Q'$.}\\
\iten  We first prove that $Q'$ is a Dyck path that stays above $P$. First note that $\delta_k(P,Q')=\delta_k(P,Q)-1$ if $i\leq k<j$ and  $\delta_k(P,Q')=\delta_k(P,Q)$ otherwise. If $\delta_k(P,Q')<0$, then $i\leq k<j$, hence $u_k$ is a descendant of $u_i$. Since the value of $\delta_r(P,Q)$ is weakly increasing along the branches of $T$, we have $\delta_k(P,Q)\geq \delta_i(P,Q)=\delta>0$, hence  $\delta_k(P,Q')\geq 0$. Thus for all $k=0,\ldots,n$,  $\delta_k(P,Q')\geq 0$, that is, $Q'$ stays above $P$.\\
\iten  We now prove that $P\lT Q'$. Suppose there exist $k,l$, such that $\delta_k(P,Q')> \delta_l(P,Q')$ with $u_k$ parent of $u_l$. Since $\delta_k(P,Q)\leq \delta_l(P,Q)$,  we have $k<i\leq l<j$. Since a vertex $u_r$ is a descendant of $u_i$ if and only if $i< r<j$, the only possibility is $l=i$. Moreover, since $u_k$ is the parent of $u_i$ we have $e_k(P)<e_i(P)=e$, hence by the choice of $e$, $\delta_k(P,Q)<\delta=\delta_i(P,Q)$. Hence, $\delta_k(P,Q')= \delta_k(P,Q)\leq \delta_i(P,Q)-1=\delta_i(P,Q')$.  We reach a contradiction. Thus $\delta_k(P,Q')\leq \delta_l(P,Q')$  whenever $u_k$ is the parent of $u_l$. By the induction hypothesis, this implies $P\lT Q'$.\\
\iten It remains to prove that  $Q'$ is covered by $Q$ in the Tamari lattice. It suffices to prove that $i\underder{Q'}j$.
Recall that for all $r=0,\ldots,n-1$, $e_r(P)+1$ is the height of the vertex $u_r$ in the tree $T$. For all $i<r<j$, the vertex $u_r$ is a descendant of $u_i$ hence $e_r(P)>e_i(P)$. Moreover, since the value of $\delta_x(P,Q)$ is weakly increasing along the branches of $T$, $\delta_r(P,Q)\geq \delta_i(P,Q)$ for all $i<r<j$. Thus, for all $i<r<j$, $e_r(Q)=e_r(P)+\delta_r(P,Q)>e_i(Q)=e_i(P)+\delta_i(P,Q)$ and $e_r(Q')=e_r(Q)-1>e_i(Q')=e_i(Q)-1$. It only remains to show that $e_j(Q')\leq e_i(Q')$. The vertex $u_j$ is the first vertex not descendant of $u_i$ around $T$, hence  $e_j(P)\leq e_i(P)$. Moreover  $\delta_j(P,Q)\leq \delta =\delta_i(P,Q)$. Furthermore, the equalities $e_i(P)= e_j(P)$ and $\delta_j(P)=\delta$ cannot hold simultaneously by the choice of $i$. Thus, $e_j(Q)=e_j(P)+\delta_j(P,Q)<e_i(Q)=e_i(P)+\delta_i(P,Q)$ and $e_j(Q')=e_j(Q)\leq e_i(Q')=e_i(Q)-1$. 
\findem

\section{Intervals of the Kreweras lattice}\label{section:Kreweras}
In this section, we study the restriction of the bijection $\Phi$ to the Kreweras intervals. 

\begin{thm}\label{thm:kreweras}
The mapping $\Phi$ induces a bijection between the intervals of the Kreweras lattice $\LK_n$ and realizers of size $n$ which are both minimal and maximal. 
\end{thm}

\begin{figure}[ht]
\begin{center}
\input{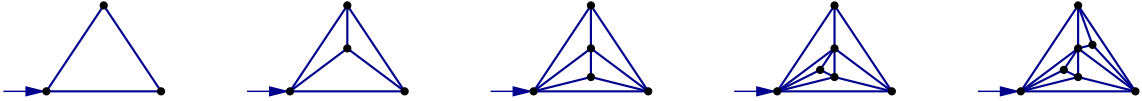}
\caption{A stack triangulation is obtained by recursively inserting a vertex of degree 3.}\label{fig:stack-triangulation}
\end{center}
\end{figure}

Before commenting on Theorem \ref{thm:kreweras}, we characterize the realizers which are both minimal and maximal. Recall that a triangulation is \emph{stack} if it is obtained from the map reduced to a triangle by recursively inserting a vertex of degree 3 in one of the (triangular) internal face. An example is given in Figure \ref{fig:stack-triangulation}.

\begin{prop}\label{prop:minandmax-is-stack} 
A realizer $R$ is both minimal and maximal if and only if the underlying triangulation $M$ is stack. (In this case, $R$ is the unique realizer of $M$.)
\end{prop}

The proof of Proposition~\ref{prop:minandmax-is-stack} uses the following Lemma.
\begin{lemma}\label{lem:restriction-realizer}
Let $M$ be a triangulation and let $R=(T_0,T_1,T_2)$ be one of its realizers. Suppose that $M$ has an internal vertex $v$ of degree 3 and let $M'$ be obtained from $M$ by removing $v$ (and the incident edges). Then, the restriction of the realizer $R$ to the triangulation $M'$ is a realizer.
\end{lemma}

\dem
By Schnyder condition, the vertex $v$ is incident to three tails and no head, hence it is a leaf in each of the trees $T_1,~T_2,~T_3$. Thus, the \emph{tree condition} is preserved by the deletion of $v$. Moreover, deleting $v$ does not deprive any other vertex of an $i$-tail, hence the \emph{Schnyder condition} is preserved by the deletion of $v$.\findembis{0cm}

\titre{Proof of Proposition \ref{prop:minandmax-is-stack}:}   \\
\ite We first prove that any realizer $R$ of a stack triangulation $M$ is minimal and maximal, that is, contains neither a cw- nor a ccw-triangle. We proceed by induction on the size of $M$. If $M$ is reduced to the triangle, the property is obvious. Let $M$ be a stack triangulation not reduced to the triangle. By definition, the triangulation $M$ contains an internal vertex $v$ of degree 3 such that the triangulation $M'$ obtained from $M$ by removing $v$ is stack. By Lemma \ref{lem:restriction-realizer}, the restriction of the realizer $R$ to $M'$ is a realizer. Hence, by the induction hypothesis, the triangulation $M'$ contains neither a cw- nor a ccw-triangle. Thus, if $C$ is either a cw- or a ccw-triangle of $M$, then $C$ contains $v$. But this is impossible since $v$ is incident to no head. \\
\ite We now prove that any realizer $R$ of a non-stack triangulation $M$ contains either a  cw- or a ccw-triangle.  \\
\iten We first prove that the property holds if $M$ has no internal vertex of degree 3 nor \emph{separating triangle} (a triangle which is not a face).  It is known that if $R$ contains a directed cycle, then it contains either a cw- or ccw-triangle (proof omitted; see \cite{Mendez:these}). Thus it suffices to prove that  $R$ contains a directed cycle. Let $u$ be the  third vertex of the internal triangle incident to the edge $(v_1,v_2)$. The vertex $u$ is such that $\Pa_1(u)=v_1$ and $\Pa_2(u)=v_2$ (see Figure~\ref{fig:config_3degenerebis}). The vertex $u$ has degree at least $4$ and is not adjacent to $v_0$ (otherwise one of the triangles $(v_0,v_1,u)$ or $(v_0,v_2,u)$ contains some vertices, hence is separating).  Thus, $u'=\Pa_0(u) \neq v_0$. Moreover, either $\Pa_1(u') \neq v_1$ or $\Pa_2(u')\neq v_2$, otherwise the triangle $(v_1,v_2,u')$ is separating. Let us assume that $u''=\Pa_1(u')\neq v_1$ (the other case is symmetrical). By Schnyder condition, the vertex $u''$ lies inside the cycle $C$ made of the edges $(v_0,v_1)$, $(v_1,u)$ and the 0-path from $u$ to $v_0$. By Schnyder condition, the 1-path from $u''$ to $v_1$ stays strictly inside $C$. Let  $C'$  be the cycle made of the edges $(v_1,u)$, $(u,u')$ and the 1-path from $u'$ to $v_1$. By Schnyder condition, the 2-path from $u''$ to $v_2$ starts inside the cycle $C'$, hence cut this cycle. Let $v$ be the first vertex of $C'$ on the 2-path from $u''$ to $v_2$. The vertex $v$ is incident to a 2-head lying inside $C'$, hence by Schnyder condition $v=u$. Thus, the cycle made of the edges $(u,u')$, $(u',u'')$ and the 2-path from $u''$ to $u$ is directed.
\begin{figure}[htb!]
\begin{center}
\input{config_3degenereebis.pstex_t}
\caption{The vertices $u$, $u'=\Pa_0(u)$ and $u''=\Pa_1(u')$.}\label{fig:config_3degenerebis}
\end{center}
\end{figure}

\iten We now prove that the property holds for any non-stack triangulation $M$ without internal vertex of degree 3. If $M$ has no separating triangle then, by the preceding point, the realizer $R$ contains either a cw- or ccw-triangle. Suppose now that $M$ has a separating triangle $\Delta$. We can choose $\Delta$ not containing any other separating triangle.  In this case, the triangulation $M'$ lying inside the triangle $\Delta$ has no separating triangle and is not stack (since no internal vertex has degree 3). Let $t_0,t_1,t_2$  be the vertices of the triangle $\Delta$. By definition, there are some vertices lying inside the triangle $\Delta$. By Lemma \ref{lem:pigeon-hole}, there is no tail incident to $\Delta$ and lying inside $\Delta$. Thus, for $i=1,2,3$, the half-edges incident to the vertex $t_i$ and lying inside $\Delta$ are heads. Moreover, the Schnyder condition implies that all the heads incident to $t_i$ have the same color. Furthermore, for each color $i=1,2,3$ there is an $i$-head incident to one of the vertices $t_0,t_1,t_2$, otherwise the vertices inside $\Delta$ would not be connected to $v_i$ by an $i$-path. Hence, we can assume without loss of generality that for $i=1,2,3$, the heads incident to $t_i$ and lying inside $\Delta$ are of color~$i$. Thus, the restriction $R'$ of $R$ to the triangulation lying inside $\Delta$ is a realizer. By the preceding point, the realizer $R'$ contains either a cw- or ccw-triangle, hence so do~$R$. \\
\iten We now prove that the property holds for any non-stack triangulation $M$.  Let $R$ be a realizer of a non-stack triangulation $M$. 
Let $M'$ be the triangulation obtained from $M'$  by recursively deleting every internal vertex of degree 3. The triangulation $M'$ is not stack and has no internal vertex of degree 3. Moreover, by Lemma \ref{lem:restriction-realizer}, the restriction $R'$ of the realizer $R$ to the triangulation $M'$ is a realizer. By the preceding point, the realizer $R'$ contains either a cw- or ccw-triangle, hence so do~$R$.
\findem

\begin{figure}[htb!]
\begin{center}
\input{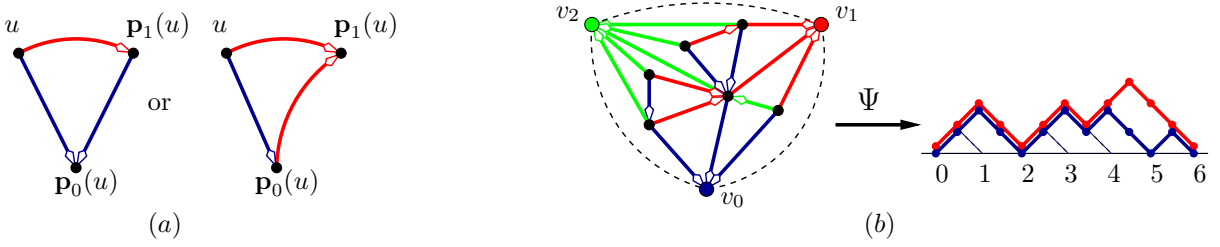}
\caption{(a) Condition for a realizer to be both minimal and maximal: $\Pa_0(\Pa_1(u))\!=\!\Pa_0(u)$ or $\Pa_1(\Pa_0(u))\!=\!\Pa_1(u)$. (b) A minimal and maximal realizer and its image by $\Psi$.}\label{fig:minimaximality}
\end{center}
\end{figure}

Given Theorem \ref{thm:kreweras} and Proposition \ref{prop:minandmax-is-stack}, the mapping $\Phi$ induces a bijection between the intervals of the Kreweras lattice and the stack triangulations. Stack triangulations are known to be in bijection with ternary trees (see for instance \cite{Zickfeld-Ziegler:poster-stack-triang}), hence we obtain a new proof that the number of intervals in  $\LK_n$ is $\frac{1}{2n+1}{3n \choose n}$. The rest of this section is devoted to the proof of Theorem \ref{thm:kreweras}. We first recall a characterization of the realizers which are both minimal and maximal. This characterization which is illustrated in Figure~\ref{fig:minimaximality} follows immediately from the characterizations of minimality and of maximality given in \cite{Bonichon:BGH02}. 

\begin{prop}[\cite{Bonichon:BGH02}] \label{prop:branch2}
A realizer $R=(T_0,T_1,T_2)$ is both minimal and maximal if and only if for any internal vertex~$u$, either $\Pa_0(\Pa_1(u))=\Pa_0(u)$ or $\Pa_1(\Pa_0(u))=\Pa_1(u)$.
\end{prop}

\begin{figure}[ht]
\begin{center}
\input{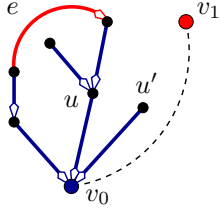}
\caption{A 1-obstruction between the vertices $u$ and $u'$.}\label{fig:obstruction}
\end{center}
\end{figure}

Let $R=(T_0,T_1,T_2)$ be a realizer of a triangulation $M$ and let $u$, $u'$ be two vertices distinct from $v_0$ and~$v_2$. We say that there is a \emph{1-obstruction} between $u$ and $u'$ if there is a 1-edge $e$ such that  the tail of $e$ appears before the first corner of $u$ while its head appears strictly between the first corner of $u$ and the first corner of $u'$ around the tree  $\bbar{T_0}$.
This situation is represented in Figure \ref{fig:obstruction}. Using Proposition \ref{prop:branch2}, we obtain the following property satisfied by realizers which are both minimal and maximal.

\begin{lemma}\label{lem:Qrelation-minmax}
Let $R=(T_0,T_1,T_2)$ be a minimal and maximal realizer and let $(P,Q)=\Psi(R)$. Let \\$v_0,u_0,u_1,\ldots,u_n\!=\!v_1$ be the vertices of the tree $\bbar{T_0}$ in clockwise order. Then, for all indices $0\leq i<j \leq n$, the relation $i\under{Q} j$ holds  if and only if the three following properties are satisfied:\\
\indent (1) the vertex $u_j$ is an ancestor of $u_i$ in the tree $T_1$,  \\
\indent (2) either $\Pa_1(\Pa_0(u_i))=u_j$ or $\Pa_0(u_i)=\Pa_0(u_j)$ (with the convention that  $\Pa_0(u_n)=v_0$),\\
\indent (3) there is no 1-obstruction between $u_i$ and $u_j$.\\
\end{lemma}

The proof Lemma \ref{lem:Qrelation-minmax} of is based on the following result.

\begin{lemma}\label{lem:Qstrongrelation-minmax}
Let $R=(T_0,T_1,T_2)$ be a minimal and maximal realizer and let $(P,Q)=\Psi(R)$. Let \\ $v_0,u_0,u_1,\ldots,u_n\!=\!v_1$ be the vertices of the tree $\bbar{T_0}$ in clockwise order. For all indices $0\leq i<j\leq n$, the relation $i\underder{Q}j$ holds if and only if $\Pa_1(u_i)=u_j$. Moreover, in this case $e_j(Q)=e_i(Q)$ if and only if $\Pa_0(u_i)=\Pa_0(u_j)$ and there is no 1-edge whose head is incident to $u_j$ and whose tail appears before the first corner of $u_i$.
\end{lemma}

\begin{figure}[h!]
\begin{center}
\input{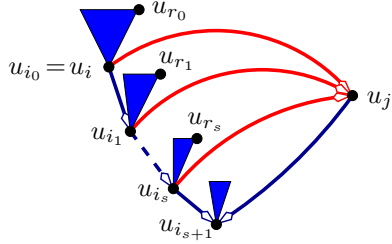}
\caption{Notations for the proof of Lemma \ref{lem:Qstrongrelation-minmax}.}\label{fig:Qstrongrelation-minmax}
\end{center}
\end{figure}

\dem Let $0\leq i<j\leq n$ such that $\Pa_1(u_i)=u_j$. \\
\ite We first prove that \emph{for all index $r=i+1,\ldots,j-1$, the inequality $e_k(Q)>e_i(Q)$ holds}.\\ 
Let $u_{i_1},\ldots,u_{i_s}$ be the vertices on the 0-path from $u_{i_0}=u_i$ to $u_{i_{s+1}}=\Pa_0(u_j)$ (that is, $\Pa_0(u_{i_k})=u_{i_{k+1}}$ for all $k=0,\ldots,s$); see Figure \ref{fig:Qstrongrelation-minmax}. The characterization of minimal and maximal realizers given in Proposition~\ref{prop:branch2} implies that $\Pa_1(u_{i_k})=u_j$ for all $k=1,\ldots,s$. For all $k=0,\ldots,s$, we denote by $r_k$ the index of the last descendant of $u_{i_k}$ around $\bbar{T_0}$ and we denote $r_{s+1}=j-1$.  Note that, for all $k=0,\ldots,s$, the vertices $u_{r_k+1},\ldots,u_{r_{k+1}}$ are descendants of $u_{i_{k+1}}$ in $T_0$. Hence, for all $k=0,\ldots,s$ and all $l=r_k+1,\ldots,r_{k+1}$, the inequality $e_l(P)> e_i(P)-k$ holds (since for any index $h$ the value $e_h(P)+1$ is the height of the vertex $u_h$ is the tree $T_0$). By the minimality condition given by Proposition \ref{prop:branch}, none of the 1-tails available at $u_i$  is matched to one the vertices  $u_{i+1},\ldots,u_{r_0}$ (since these vertices are descendants of $u_i$). Moreover, none of these available 1-tails is matched to one of the  vertices  $u_{r_0+1},\ldots,u_{j-1}$ or there would be a crossing with the 1-edge $(u_i,u_j)$  (see Figure \ref{fig:Qstrongrelation-minmax}). Hence, the 1-tails available at $i$ are also available at all the vertices $u_{i+1},\ldots,u_{j}-1$. Moreover, for all $k=1,\ldots,s$ and all $l=r_k+1,\ldots,r_{k+1}$ the $k$ 1-tails incident to each of the vertices $u_{i_0},\ldots,u_{i_{k-1}}$ are available at the vertex $u_l$. Thus, given Lemma \ref{lem:nb-available-tails}, for all $k=0,\ldots,s$, for all $l=r_k+1,\ldots,r_{k+1}$, $\delta_l(P,Q)\leq \delta_i(P,Q)+k$. Thus, for all $l=i+1,\ldots,j-1$, the inequality $e_l(Q)=e_l(P)+\delta_l(P,Q)>e_i(Q)=e_i(P)+\delta_i(P,Q)$ holds.\\ 
\ite It only remains to prove that \emph{the inequality $e_j(Q)\leq e_i(Q)$ holds  and equality occurs if and only if $\Pa_0(u_i)=\Pa_0(u_j)$ and there is no 1-edge whose head is incident to $u_j$ and whose tail appears before the first corner of~$u_i$}.\\
\iten Since the realizer $R$ is minimal the vertex $\Pa_0(u_j)$ is an ancestor of $u_j$ in the tree $T_0$ (by Proposition \ref{prop:branch}). Hence, the inequality $e_j(P)\leq e_i(P)$ holds  and equality occurs if and only if $\Pa_0(u_i)=\Pa_0(u_j)$. We now compare the values of $\delta_i(P,Q)$ and $\delta_j(P,Q)$ which are the number of tails available at $u_i$ and at $u_j$ respectively (by Lemma \ref{lem:nb-available-tails}). \\
\iten We first prove that \emph{any 1-tails available at $u_j$ is also available at $u_i$}. No 1-tail available at $u_j$ is incident to a vertex $u_l$ with $r_0<l<j$ or the corresponding 1-edge would cross the edge $(u_i,u_j)$ (see Figure \ref{fig:Qstrongrelation-minmax}). Moreover,  the characterization of minimal and maximal realizers given in Proposition~\ref{prop:branch2} implies that no  1-tail available at $u_j$ is incident to a vertex $u_l$ with $i<l\leq r_0$ (since these vertices are descendants of $u_i$).  Hence, any the 1-tail available at $u_j$ is also available at $u_i$. \\
\iten We now prove \emph{any 1-tail available at $u_i$ is available at $u_j$ except if the corresponding 1-head is incident to $u_j$}. Clearly, no 1-tail available at $u_i$ is such that the corresponding 1-head is incident to a vertex $u_l$ with $r_0<l<j$ or the 1-edge under consideration would cross the edge $(u_i,u_j)$ (see Figure \ref{fig:Qstrongrelation-minmax}). Since the realizer $R$ is minimal, there is no 1-tail available at $u_i$ and such that the corresponding 1-head is a vertex  $u_l$ with $i<l\leq r_0$ (since these vertices are descendants of $u_i$). Hence, if a 1-tail available at $u_i$ is not available at $u_j$, then the corresponding 1-head is incident to $u_j$.\\
\iten Given Lemma \ref{lem:nb-available-tails}, the preceding points imply that the inequality $\delta_i(P,Q)\leq\delta_j(P,Q)$ holds and equality occurs if and only if there is no 1-edge whose head is incident to $u_j$ and whose tail appears before the first corner of $u_i$. Hence,  $e_j(Q)=e_j(P)+\delta_j(P,Q)\leq e_i(Q)=e_i(P,Q)+\delta_i(P,Q)$  and equality occurs if and only if $\Pa_0(u_i)=\Pa_0(u_j)$ and no index $k<i$ is such that $\Pa_1(u_k)=u_j$.  
\findem

\begin{figure}[h!]
\begin{center}
\input{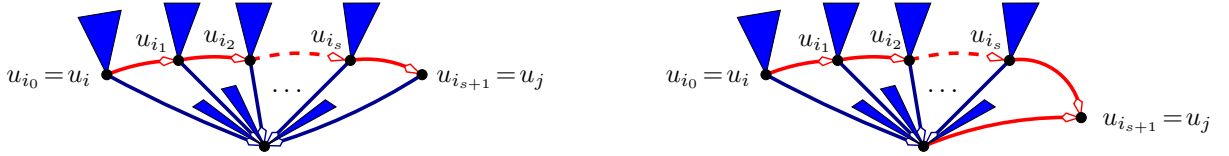}
\caption{Notations for the proof of Lemma \ref{lem:Qrelation-minmax}.}\label{fig:Qrelation-minmax}
\end{center}
\end{figure}

\demof{Lemma~\ref{lem:Qrelation-minmax}}\\
\ite We suppose  that $i\under{Q} j$ and we want to prove the properties (1), (2) and (3). \\
(1) Since  $i\under{Q} j$, there are indices $i_0=i,i_1,\ldots,i_{s+1}=j$ such that $i_0\underder{Q} i_1\underder{Q}\ldots \underder{Q} i_{s+1}$ and  $e_{i_0}(Q)=e_{i_1}(Q)=\cdots=e_{i_s}(Q)$. Lemma  \ref{lem:Qstrongrelation-minmax} implies $\Pa_1(u_{i_k})=u_{i_{k+1}}$ for all $k\leq s$. Hence $u_j$ is a ancestor of $u_i$ in the tree $T_1$. \\
(2) Since  $e_{i_k}(P)=e_i(P)$ for all $k\leq s$, Lemma  \ref{lem:Qstrongrelation-minmax} implies $\Pa_0(u_{i_k})=\Pa_0(u_i)$ for all $k\leq s$. Moreover, $\Pa_1(u_{i_s})=u_j$, thus  Proposition~\ref{prop:branch2} implies that either $\Pa_0(u_j)=\Pa_0(u_{i_s})=\Pa_0(u_i)$ or $u_j=\Pa_1(\Pa_0(u_{i_s}))=\Pa_1(\Pa_0(u_i))$.  This situation is represented in Figure \ref{fig:Qrelation-minmax}. \\
(3) We want to prove that there is no 1-obstruction between $u_i$ and $u_j$. We suppose that the tail of a 1-edge $e$ appears before the first corner of $u_i$ around $T_0$ and we want to prove that the corresponding 1-head $h$ is not incident to a vertex $u_k$ with $i<k<j$. Clearly, if the 1-head $h$ is incident to a vertex $u_k$ with $i<k<j$, then the vertex $u_k$ is either one of the vertices $u_{i_0},u_{i_2},\ldots,u_{i_s}$ or one of their descendants (otherwise, the edge $e$ would cross one of the 1-edges $(u_{i_0},u_{i_1})$, \ldots ,$(u_{i_s},u_{i_{s+1}})$; see Figure \ref{fig:Qrelation-minmax}). Since $e_{i_0}(Q)=e_{i_1}(Q)=\cdots=e_{i_s}(Q)$, Lemma  \ref{lem:Qstrongrelation-minmax} implies that $u_k$ is none of the vertices $u_{i_1},u_{i_2},\ldots,u_{i_s}$. Moreover, since the realizer $R$ is minimal, Proposition~\ref{prop:branch} implies that $u_k$ is not a (proper) descendant of one of the vertices $i_0,\ldots,i_s$. Thus, the 1-head $h$ is not incident to a vertex $u_k$ with $i<k<j$ and $e$ is not creating a 1-obstruction.\\
\ite We suppose that the vertices $u_i$ and $u_j$ satisfy the properties (1), (2) and (3) and want to prove that $i\under{Q}j$. Observe first that by property (1), there are indices $i_0=i,i_1,\ldots,i_{s+1}=j$  such that $\Pa_1(u_{i_k})=u_{i_{k+1}}$.\\
\iten We first prove that, \emph{for all $k=1,\ldots,s$, $\Pa_0(u_{i_k})=\Pa_0(u_{i})$}; this situation is represented in Figure \ref{fig:Qrelation-minmax}. \\
Suppose the contrary and consider the first index $k\in \{1,\ldots,s\}$ such that $\Pa_0(u_{i_{k}})\neq \Pa_0(u_i)$. In this case, $u_{i_k}=\Pa_1(u_{i_{k-1}})$ and  $\Pa_0(u_{i_k})\neq \Pa_0(u_{i_{k-1}})= \Pa_0(u_i)$. Since the realizer $R$ is minimal and maximal, Proposition~\ref{prop:branch2} implies that  $u_{i_{k}}=\Pa_1(\Pa_0(u_{i_{k-1}})=\Pa_1(\Pa_0(u_{i}))$. Thus, the vertices $u_{i_{k+1}},\ldots,u_{i_{s+1}}$ are distinct from $\Pa_1(\Pa_0(u_i))$ and are ancestors of $\Pa_0(u_i)$ in the tree $T_1$. In particular, $u_j=u_{i_{s+1}}\neq \Pa_1(\Pa_0(u_i))$, and $\Pa_0(u_j)\neq \Pa_0(u_i)$. This contradicts Property~(2).\\ 
\iten We now prove that \emph{for all index $k=1,\ldots,s$ there is no 1-edge $e$ whose head is incident to $u_{i_k}$ and whose tail appears before the first corner of $u_{k-1}$}. Suppose that such a 1-edge $e$ exist.
Observe that the 1-tail $t$ of the edge $e$ do not appear before the first corner of $u_i$ otherwise the edge $e$ creates a 1-obstruction between $u_i$ and $u_j$. Hence, the 1-tail $t$ is incident either to one of the vertices $u_{i_0},\ldots,u_{i_{k-2}}$ or to one of their descendants  (otherwise, the edge $e$ would cross one of the 1-edges $(u_{i_0},u_{i_1})$, \ldots ,$(u_{i_{k-2}},u_{i_{k-1}})$; see Figure \ref{fig:Qrelation-minmax}). Moreover, the  1-tail $t$ is not incident to the vertices $u_{i_0},\ldots,u_{i_{k-2}}$, otherwise $e$ would create a cycle in the tree $T_1$. Lastly, since the realizer $R$ is minimal, the 1-tail $t$ is not incident to a descendant of $u_{i_l},~l=0,\ldots,k-2$. Thus the 1-tail $t$ does not appear before the first corner of $u_{k-1}$.\\
\iten By Lemma \ref{lem:Qstrongrelation-minmax}, the preceding points imply $i_k\underder{Q} i_{k+1}$ and  $e_{i_k}(Q)=e_{i}(Q)$ for all $k=0\ldots s$. Thus, $i\under{Q}j$.\findembis{0cm}

\begin{figure}[ht]
\begin{center}
\hspace{-1cm}\input{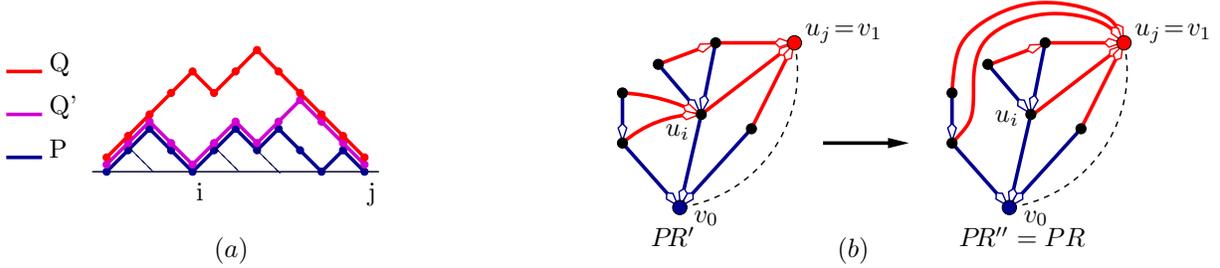}
\caption{(a) The Dyck paths $P\lK Q'\lK Q$. (b) The prerealizer $\PR''$ is obtained from $\PR'=(T_0,T_1')$ by moving $\be_i'$ 1-heads from the first corner of $u_i$ to the first corner of $u_j$.}\label{fig:transfer-1-heads}
\end{center}
\end{figure}

\demof{Theorem \ref{thm:kreweras}}
Let $P=NS^{\al_1}\ldots NS^{\al_n}$ and $Q=NS^{\be_1}\ldots NS^{\be_n}$ be two Dyck paths and let $R=(T_0,T_1,T_2)=\Phi(P,Q)$. Let $v_0,u_0,u_1,\ldots,u_n=v_1$ be the vertices of the tree $\bbar{T_0}$ in clockwise order.\\
\ite We suppose that $P\lK Q$ and we want to prove that the realizer $R$ is minimal and maximal. We proceed by induction on $\Delta(P,Q)$.\\
\iten We first suppose that $\Delta(P,Q)=0$, that is $P=Q$, and we want to prove that $R$ is minimal and maximal. 
Let $\mW$ be the word obtained by making the tour of $\bbar{T_0}$ and writing $N$ (resp. $S$) when following an edge of $\bbar{T_0}$ for the first (resp. second) time and writing $\BN$ (resp. $\BS$) when crossing a 1-tail (resp. 1-head). By definition of the mapping $\omega$, the restriction of $\mathcal{W}$ to the letters $N,S$ is $\omega(\bbar{T_0})=NS^{\al_1}\ldots NS^{\al_n}NS$. Moreover, for all $i=0,\ldots,n$ there are $\al_i$ 1-heads incident to the first corner of $u_i$ and one 1-head incident to its last corner. Thus, $\mathcal{W}=N(\BN S)^{\al_1}N\BS^{\al_1}(\BN S)^{\al_2} \ldots N\BS^{\al_{n-1}}(\BN S)^{\al_n}N\BS^{\al_{n}}S$. Between any letter $\BN$ of $\mW$ and the corresponding letter $\BS$ there is exactly one letter $N$. 
Thus, for any internal vertex $u$, the vertex $\Pa_1(u)$ is the first vertex appearing after the last corner of $u$ around $\bbar{T_0}$ (that is,  the first vertex which is not a descendant of $u$ appearing after $u$  around $\bbar{T_0}$). By Proposition \ref{prop:branch2}, this implies that $R$ is minimal and maximal. \\
\iten We now suppose that  $\Delta(P,Q)>0$.  In this case, there is a Dyck path $Q'=NS^{\be_1'}\ldots NS^{\be_n'}$ covered by $Q$ in the Kreweras lattice and such that $P\lK Q'$. Since $Q'$ is covered by $Q$ is the Kreweras lattice, there are indices $0\leq i<j \leq n$ such that $i\under{Q'}j$ and $\be_i=0$, $\beta_j=\beta_i'+\beta_j'$ and $\be_k=\be_k'$ for all $k\neq i,j$ (this situation is represented in Figure \ref{fig:transfer-1-heads}~(a)). By the induction hypothesis, the realizer $R'=(T_0',T_1',T_2')=\Phi(P,Q')$ is both minimal and maximal. Moreover, by definition of the bijection $\Phi$, the trees $T_0$ and $T_0'$ are the same. We use this fact to identify the vertices in the prerealizers $\PR=(T_0,T_1)$ and $\PR'=(T_0,T_1')$ that we denote by $v_0,u_0,u_1,\ldots,u_n=v_1$ in clockwise order around $\bbar{T_0}=\bbar{T_0'}$. We also denote by $\Pa_1'(u)$ the parent of any vertex $u$ in $T_1'$. 
\begin{itemize}
\item We first prove that \emph{for any vertex $v$,  $\Pa_1'(v)=\Pa_1(v)$ except if  $\Pa_1'(v)=u_i$ in which case $\Pa_1(v)=u_j$.}
Since $i\under{Q} j$, Lemma \ref{lem:Qrelation-minmax} implies that there is no 1-obstruction between $u_i$ and $u_j$ in the realizer~$R'$. Thus, the $\be_i'$ 1-heads incident to $u_i$ can be unglued from the first corner of $u_i$ and glued to the first corner of $u_j$ without creating any crossing in the prerealizer $\PR'=(T_0,T_1')$ (the transfer of  the $\be_i'$ 1-heads is represented in Figure \ref{fig:transfer-1-heads}~(b)). Let $\PR''=(T_0,T_1'')$ be the colored map obtained. Clearly, $\PR''=(T_0,T_1'')$ satisfies the \emph{tree condition} ($T_1''$ is a tree), the \emph{corner condition} (the 1-heads are in first corners, the 1-tails are in last corners) and the \emph{order condition} (any 1-tail appears before the corresponding 1-head around $\bbar{T_0}$), therefore $\PR''$ is a prerealizer. Moreover, for all $i=0,\ldots,n$, there are $\beta_i$ 1-heads incident to the vertex $u_i$. Thus, by definition of the mapping~$\Phi$, the prerealizer $\PR''$ is equal to $\PR=(T_0,T_1)$. Since the only difference between the prerealizers $\PR'$ and $\PR$ is that the 1-heads incident to $u_i$ in $\PR'$ are incident to $u_j$ in $\PR$, the property holds.
\item We now prove that \emph{the realizer $R=(T_0,T_1,T_2)$ is minimal and maximal.} If the realizer $R$ is not both minimal and maximal, there is a vertex $u$ such that $\Pa_1(u)\neq \Pa_1(\Pa_0(u))$ and $\Pa_0(\Pa_1(u))\neq\Pa_0(u)$. Since the realizer $R'$ is both minimal and maximal, either $\Pa_1'(u)=\Pa_1'(\Pa_0(u))$ or $\Pa_0(\Pa_1'(u))=\Pa_0(u)$. But  $\Pa_1'(u)\neq \Pa_1'(\Pa_0(u))$, otherwise $\Pa_1(u)=\Pa_1(\Pa_0(u))$. Thus, $\Pa_0(\Pa_1'(u))=\Pa_0(u)$ and $\Pa_1'(u)=u_i$. Hence, $\Pa_0(u_i)=\Pa_0(u)$ and $\Pa_1(u)=u_j$.  Moreover, since $i\under{Q'} j$,  Lemma~\ref{lem:Qrelation-minmax} implies that either $\Pa_0(u_i)=\Pa_0(u_j)$ or $\Pa_1'(\Pa_0(u_i))=u_j$. But, if $\Pa_0(u_i)=\Pa_0(u_j)$,  then $\Pa_0(u)=\Pa_0(u_i)=\Pa_0(u_j)=\Pa_0(\Pa_1(u))$ which is forbidden. And, if  $\Pa_1'(\Pa_0(u_i))=u_j$, then $\Pa_1(\Pa_0(u))=\Pa_1(\Pa_0(u_i))=\Pa_1'(\Pa_0(u_i))=u_j=\Pa_1(u)$  which is also forbidden. We reach a contradiction.
\end{itemize}

\begin{figure}[ht]
\begin{center}
\input{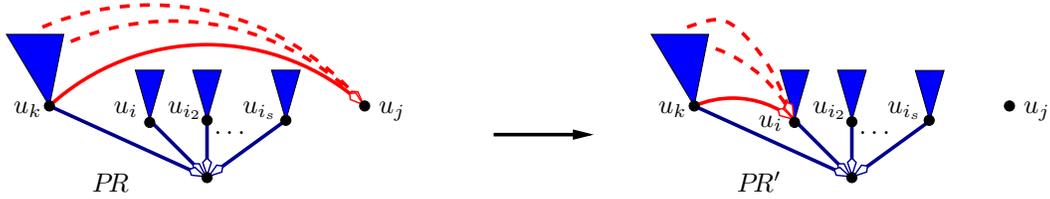}
\caption{The vertices $u_k,u_i,u_j$ in the prerealizer $\PR=(T_0,T_1)$ and $\PR'=(T_0,T_1')$.}\label{fig:transfer-1-heads-back}
\end{center}
\end{figure}

\ite We suppose that the realizer $R$ is minimal and maximal and we want to prove that $P\lK Q$. We proceed by induction on $\Delta(P,Q)$. If $\Delta(P,Q)=0$, then $P=Q$ and the property holds. We suppose now that $\Delta(P,Q)>0$ and we denote by $v_0,u_0,u_1,\ldots,u_n=v_1$  the vertices of the tree $\bbar{T_0}$ in clockwise order.\\
\iten We first prove that \emph{there are indices $0\leq k<i<j\leq n$ such that $\Pa_0(u_k)=\Pa_0(u_i)$ and $\Pa_1(u_k)=u_j$.} We suppose that no such indices exist and we want to prove that $P=Q$. 
Let $u$ be an internal vertex. If $u$ has a sibling in $\bbar{T_0}$ appearing after $u$ around $\bbar{T_0}$, then $\Pa_1(u)$ is the first such sibling (since the indices $i,j,k$ do not exist), else $\Pa_1(u)=\Pa_1(\Pa_0(u))$  (since the realizer $R$ is minimal and maximal).
Thus, for any vertex $u$, $\Pa_1(u)$ is the first vertex appearing after the last corner of $u$ around $\bbar{T_0}$. Let $\mW$ be the word obtained by making the tour of $\bbar{T_0}$ and writing $N$ (resp. $S$) when following an edge of $\bbar{T_0}$ for the first (resp. second) time and writing $\BN$ (resp. $\BS$) when crossing a 1-tail (resp. 1-head). By definition of the mapping $\omega$, the restriction of $\mathcal{W}$ to the letters $N,S$ is $\omega(\bbar{T_0})=NS^{\al_1}\ldots NS^{\al_n}NS$. Moreover, for all $i=0,\ldots,n$ there are $\be_i$ 1-heads in the first corner of $u_i$ and one 1-head in its last corner. Thus, $\mathcal{W}=N(\BN S)^{\al_1}N\BS^{\be_1}(\BN S)^{\al_2} \ldots N\BS^{\be_{n-1}}(\BN S)^{\al_n}N\BS^{\be_{n}}S$. Moreover, between any letter $\BN$ of $\mW$ and the corresponding letter $\BS$ there is exactly one letter $N$. Thus, $\be_1=\al_1$,\ldots,  $\be_n=\al_n$, that is, $P=Q$.\\
\iten Let $k<i<j$ be as described in the preceding point with $k$ maximal and $i$ minimal with respect to $k$ (i.e. $u_i$ is the first sibling of $u_k$ appearing after $u_k$ around the tree $\bbar{T_0}$). This situation is represented in  Figure~\ref{fig:transfer-1-heads-back}. Observe that no 1-head is incident to $u_i$ in the prerealizer $\PR=(T_0,T_1)$ (see Figure \ref{fig:transfer-1-heads-back}), hence $\be_i=0$,. Let $H$ be the set of 1-heads incident to $u_j$ and such that the corresponding 1-tail is either incident to $u_k$ or to one of its descendants. One can unglue the  1-heads in $H$ from the first corner of $u_j$ and glue them to the first corner of $u_i$ without creating any crossing (see Figure \ref{fig:transfer-1-heads-back}).  Moreover, the resulting colored map $\PR'$ is easily seen to be a prerealizer that we denote by $\PR'=(T_0,T_1')$. Let $R'$ be the realizer corresponding to the prerealizer $\PR'$ and let $Q'=NS^{\be_1'}\ldots NS^{\be_n'}$ be the Dyck path such that $\Phi(P,Q')=R'$. By definition of $\Phi$, we have  $\be_i'=|H|$, $\be_j'=\be_j-|H|$ and $\be_l'=\be_l$ for all $l\neq i,j$. \\
\iten We now prove that  \emph{the realizer $R'=\Phi(P,Q')$ is minimal and maximal.} By Proposition \ref{prop:branch2}, we only need to prove that for every internal vertex $u$, either  $\Pa_0(\Pa_1'(u))=\Pa_0(u)$ or $\Pa_1'(\Pa_0(u))=\Pa_1'(u)$, where $\Pa_1'(u)$ denotes the parent of $u$ in the tree $T_1'$. Suppose that there is a vertex $u$ not satisfying this condition. Note first that $u\neq u_k$ since  $\Pa_0(\Pa_1'(u_k))=\Pa_0(u_k)$. Since the realizer $R$ is minimal and maximal, either $\Pa_0(\Pa_1(u))=\Pa_0(u)$ or $\Pa_1(\Pa_0(u))=\Pa_1(u)$. Suppose first $\Pa_0(\Pa_1(u))=\Pa_0(u)$. 
In this case, the vertex $u$ is a descendant of $u_k$ (otherwise, $\Pa_0(\Pa_1'(u))=\Pa_0(\Pa_1(u))=\Pa_0(u)$), and $\Pa_1'(u)=u_j$ (for the same reason). Therefore,  $\Pa_0(u_j)=\Pa_0(\Pa_1(u))=\Pa_0(u)$ implies that $u_j$ is a descendant of $u_k$. This is impossible since $u_j$ appears after $u_i$ around $\bbar{T_0}$. Suppose now that  $\Pa_1(\Pa_0(u))=\Pa_1(u)$.  In this case, the vertex $u$ is a descendant of $u_k$ (otherwise, $\Pa_1'(\Pa_0(u))=\Pa_1(\Pa_0(u))=\Pa_1(u)=\Pa_1'(u)$), and $\Pa_1(\Pa_0(u))=\Pa_1(u)=u_j$ (for the same reason). Thus $\Pa_1'(\Pa_0(u))=\Pa_1'(u)=u_i$.  We reach again a contradiction. \\
\iten We now prove that \emph{the Dyck path $Q'$ is covered by $Q$ in the Kreweras lattice.}
By definition of the covering relation in the Kreweras lattice $\LK$, it suffices to prove that $i\under{Q'}j$. Since the realizer $R'$ is minimal and maximal, it suffices to prove that the conditions (1), (2)  and (3) of Lemma \ref{lem:Qrelation-minmax} hold. Clearly, there is no 1-obstruction between the vertices $u_i$ and $u_j$ in the realizer $R'$ (see Figure \ref{fig:transfer-1-heads-back}), hence condition (3) holds. Moreover, since  the realizer $R$ is minimal and maximal, either $\Pa_0(u_k)=\Pa_0(u_j)$ or $\Pa_1(\Pa_0(u_k))=u_j$. Thus, either  $\Pa_0(u_i)=\Pa_0(u_j)$ or $\Pa_1(\Pa_0(u_i))=u_j$, hence condition (2) holds. 
Let $i=i_1, i_2,\ldots,i_s$ be the indices of the siblings of $u_k$ appearing between $u_k$ and $u_j$ in clockwise order around $\bbar{T_0}$ (see Figure \ref{fig:transfer-1-heads-back}). By the choice of $k$, we get $\Pa_1(u_{i_r})=u_{i_{r+1}}$ for all $r<s$. Moreover, since  the realizer $R$ is minimal and maximal, either $\Pa_0(u_k)=\Pa_0(u_j)$ or $\Pa_1(\Pa_0(u_k))=u_j$. If either case, we get   $\Pa_1(u_s)=u_j$. Thus, $\Pa_1'(u_{i_r})=\Pa_1(u_{i_r})=u_{i_{r+1}}$ for all $r<s$, and $\Pa_1'(u_s)=\Pa_1(u_s)=u_j$. Hence, $u_j$ is an ancestor of $u_i$ in the tree $T_1'$, that is, condition (1) holds.\\
\iten The realizer  $R'=\Phi(P,Q')$ is minimal and maximal, hence by the induction hypothesis $P\lK Q'$.  Moreover, the path $Q'$ is covered by $Q$ in the Kreweras lattice. Thus, $P\lK Q$.
\findem

\newpage

\titre{Acknowledgments:} The authors are very grateful to Mireille Bousquet-Mélou for pointing them to the enumerative result of Frédéric Chapoton and to Xavier Viennot for fruitful discussions and suggestions.

\bibliography{../../../biblio/allref.bib}

\bibliographystyle{plain}

\end{document}